\documentclass{article}

\usepackage{amssymb}
\usepackage{amsfonts}
\usepackage{amsmath}
\usepackage{url}
\usepackage{latexsym}

\def\deg{\mathop{\rm deg }\nolimits}
\def\rank{\mathop{\rm rank}\nolimits}
\def\lcm{\mathop{\rm lcm }\nolimits}

\newcommand{\efe}{\mathbb F}
\newcommand{\FF}{\mathbb F}
\newcommand{\CC}{\mathbb C}
\newcommand{\RR}{\mathbb R}

\newcommand{\ba}{\mathbf a}
\newcommand{\bd}{\mathbf d}
\newcommand{\bg}{\mathbf g}
\newcommand{\bc}{\mathbf c}
\newcommand{\bb}{\mathbf b}
\newcommand{\bu}{\mathbf u}
\newcommand{\bv}{\mathbf v}

\newtheorem{theorem}{Theorem}[section]
\newtheorem{corollary}[theorem]{Corollary}
\newtheorem{lemma}[theorem]{Lemma}
\newtheorem{definition}[theorem]{Definition}
\newtheorem{problem}[theorem]{Problem}

\newtheorem{remark}[theorem]{Remark}
\newtheorem{example}[theorem]{Example}

\title{Row or column completion of polynomial matrices of given degree II
  \thanks{
This work was 
supported by grant PID2021-124827NB-I00 funded by MCIN/AEI/ 10.13039/501100011033 and by ``ERDF A way of making Europe'' by the ``European Union''.
The first and third authors were also supported 
by grant GIU21/020 funded by UPV/EHU.
  }
  }

\author{Agurtzane Amparan\thanks{Departamento de Matem\'aticas, Universidad del Pa\'is Vasco UPV/EHU, Bilbao, Spain, {agurtzane.amparan@ehu.eus}, {silvia.marcaida@ehu.eus}}
\and Itziar Baraga\~na\thanks{Departamento de Ciencia de la Computaci\'on e I.A., 
Universidad del Pa\'{\i}s Vasco, UPV/EHU, Donostia-San Sebasti\'an, 
Spain, {itziar.baragana@ehu.eus}}
\and Silvia Marcaida\footnotemark[2]
\and Alicia Roca\thanks{Departamento de Matem\'atica Aplicada, IMM, Universitat Polit\`ecnica de Val\`encia, 46022 Valencia, Spain,   {aroca@mat.upv.es}}}

\date{}

\begin{document}

\maketitle

\begin{abstract}
The row (column) completion problem  of polynomial matrices of given degree with prescribed eigenstructure has been studied in \cite{AmBaMaRo23}, where several results of prescription of some of the four types of invariants  that form the eigenstructure have also been obtained. In this paper we conclude the study, 
solving the completion for  the cases not covered there.
More precisely,  we solve  the row completion problem of a polynomial matrix  when we prescribe  the infinite (finite) structure and column and/or row minimal indices, and finally  
the column and/or row  minimal indices.
 The necessity of the characterizations  obtained holds to be true  over arbitrary fields in all cases, whilst to prove the sufficiency it is required, in some of the cases, to work  over algebraically closed fields.
\end{abstract}


{\bf AMS:}
 15A18, 15A54, 15A83

\section{Introduction}\label{secintro}
Numerous problems in  different scientific areas  lead  to the study  of systems of differential or difference equations, where polynomial matrices  play often an important role \cite{Bhms11, Fo75, GLR85, Kail80, Rose70, TiMe01, Vard91}. 
On the other hand, an important problem we are interested in   is the matrix completion problem, which appears for instance in system design (see for example \cite{BoDo94} and \cite{LoMoZaZaLAA98}; other references on matrix completions are analyzed below). 

These ideas led us to study in \cite{AmBaMaRo23} a polynomial matrix completion problem consisting in finding necessary and sufficient conditions for the existence of a polynomial matrix of given degree when the eigenstructure, or part of it, and some of its rows are prescribed.
By transposition, the solution to the column completion problem is also obtained. 
If the polynomial matrix is of degree one, i.e., for matrix pencils, the problem has already been solved in 
\cite{Do08} when the whole structure is prescribed (see \cite{Do13} or \cite{DoSt19} for a simplified explicit solution).

 There are many antecedents of this problem, and the literature about matrix completion results is very vast.
For polynomial matrix completion see, for example, \cite{Sa79, Th79}, where conditions of degree are not imposed. In the case of constant matrices (polynomial matrices of degree 0) see \cite{ZaLAA87} for row completion up to a square matrix with prescribed similarity invariants,  and \cite{Do05} for row completion  up to a  rectangular matrix with prescribed feedback invariants.

The eigenstructure of a polynomial matrix is formed by the finite structure (invariant factors), infinite structure (infinite elementary divisors or, equivalently, partial multiplicities of
$\infty$) and the singular structure (right and left minimal indices). 
In \cite{AmBaMaRo23}, with the aim of simplifying the expressions, instead of working with the invariant factors and  the infinite elementary divisors, we used the homogeneous invariant factors, which are uniquely determined by the finite and infinite structures. Moreover, in that paper, the right and left minimal indices were decreasingly ordered and were called column and row minimal indices, respectively.

In \cite{AmBaMaRo23}  the row completion problem of a polynomial matrix of given degree is also solved when instead of prescribing the whole  eigenstructure only part of it is prescribed, more precisely, when  we prescribe the finite and infinite structures and the
column (row) minimal indices, and also the finite and/or infinite
structures. 

In this paper, using mainly combinatorial techniques, we solve 
the
remaining cases in Section 3: 
prescription of infinite (finite) and singular structures in Theorem \ref{theoprescrdeirmicmi} (Theorem \ref{theoprescrfifrmicmi}),
prescription of infinite (finite) structure and column minimal indices in Corollary \ref{corprescrdeicmi} (Corollary \ref{corprescrfifcmi}), 
prescription of infinite (finite) structure and row minimal indices in Corollary \ref{corprescrdeirmi} (Corollary \ref{corprescrfifrmi}), 
prescription of singular structure in  Theorem \ref{theoprescrsing}, and 
prescription of row (column) minimal indices in  Corollary \ref{corprescrrmi}
(Corollary \ref{propprescrcmi}).

\section{Preliminaries}\label{secprel}

This paper completes the work started in  \cite{AmBaMaRo23}. For basic definitions and properties we refer to Section 2 of that reference.

Let $\FF$ be a field, $\FF[s]$  the ring of polynomials in the indeterminate $s$ with coefficients in $\FF$, $\FF(s)$ the field of fractions of
$\FF[s]$, and  $\FF[s, t]$   the ring of polynomials in two
variables $s, t$, with coefficients in $\FF$.
A polynomial in $\FF[s]$ is \textit{monic} if its leading coefficient is 1.  We say that a polynomial in $\FF[s,t]$ is {\em monic} if it is monic with respect to the variable $s$. Given two polynomials $p, q$, by $p\mid q$ we mean that $p$ is a divisor of $q$, and  $\lcm(p,q)$ is the monic least common multiple of $p$ and $q$. 
   
As usual, let  $\FF^{m\times n}$ and $\FF[s]^{m\times n}$ be the vector spaces over $\FF$  of $m \times n$ matrices with elements in $\FF$ and  $\FF[s]$, respectively.

Let  $P(s)\in \FF[s]^{m\times n}$ be a polynomial matrix. It can be written as $P(s)=P_ds^d+P_{d-1}s^{d-1}+\cdots+P_1s+P_0\in\efe[s]^{m\times n}$, $P_i\in  \FF^{m\times n}$,  $0\leq i \leq  d$, with $P_d\neq 0$,  for some integer $d$. Then, $d$ is  the {\em degree} of $P(s)$ (denoted by $\deg(P(s))$). 

The {\em normal rank} of $P(s)$, denoted by $\rank (P(s))$,  is the order of the largest non identically zero minor of $P(s)$, i.e., it is the rank of $P(s)$ considered as a matrix on  $\FF(s)$. We will refer to it as the {\em rank} of $P(s)$.

As mentioned in the introduction, the {\em eigenstructure} of a polynomial matrix $P(s)\in  \FF[s]^{m\times n}$, $\rank(P(s))=r$, is formed by the  {\em finite structure}, i.e., a chain of monic polynomials
$\alpha_1(s)\mid \dots \mid \alpha_r(s)$ called {\em invariant factors}, the  {\em infinite structure}, i.e., a sequence of non negative integers $e_1\leq \dots \leq e_r$ called {\em partial multiplicities of $\infty$}, and the {\em singular structure}, i.e., two sequences of non negative integers $c_1\geq\dots \geq c_{n-r}$ and  $u_1\geq\dots \geq u_{m-r}$ called {\em column} and {\em row minimal indices}, respectively.
In the literature about polynomial matrices the column and row minimal  indices  are usually increasingly ordered  and they  are called {\em right} and {\em left minimal indices}, respectively. 
 If two  polynomial matrices are strictly equivalent, then they have the same eigenstructure,  but, unlike in matrix pencils, the converse is not true in general (see \cite{DeDoMa14}).

The invariant factors and the  partial multiplicities of $\infty$ can be  summarized in a chain of monic homogeneous polynomials $\phi_1(s,t) \mid \dots \mid \phi_r(s,t)$, $\phi_i(s,t) \in \FF[s,t]$, $1\leq i\leq r$, called the  {\em homogeneous invariant  factors} of $P(s)$, of the form 
\begin{equation}\label{hif}
\phi_i(s,t)=t^{e_i}t^{\deg(\alpha_i)}\alpha_i\left(\frac st\right), \quad 1\leq i\leq r.
\end{equation}
In turn, the homogeneous invariant factors of a polynomial matrix determine the invariant factors $\alpha_i(s)=\phi_i(s,1),$  $1\leq i \leq r$,
and the partial multiplicities of $\infty,$ $e_1 \leq\dots \leq e_r$.

We will take $\phi_i(s,t)=1$ whenever $i<1$ and $\phi_i(s,t)=0$  when $i>r$. As a consequence, $\alpha_i(s)=1$ and $e_i=0$ for $i<1$, and  $\alpha_i(s)=0$ for $i>r$. 
 We also agree that $e_i=+\infty$ for $i>r$.

The following result is in \cite{DeDoMa14}; we state it here  in terms of the homogeneous invariant factors.

\begin{lemma} {\rm (Index Sum Theorem for Matrix Polynomials \cite[Theorem 6.5]{DeDoMa14})}
\label{theoDeDoMa}
 Let $P(s)\in\efe[s]^{m\times n}$ be a polynomial matrix, $\deg(P(s))=d$, 
 $\rank(P(s))=r$.
Let $\phi_1(s,t)\mid\cdots\mid\phi_r(s,t)$, 
$c_1\geq\cdots\geq c_{n-r}$ and $u_1\geq\cdots\geq u_{m-r}$
be the homogeneous invariant factors, column minimal indices and row minimal indices of $P(s)$, respectively. Then, 
\begin{equation}\label{eqIST}
\sum_{i=1}^r \deg(\phi_i)+\sum_{i=1}^{n-r} c_i+\sum_{i=1}^{m-r} u_i=rd.
\end{equation}
\end{lemma}

Notice that if $P(s)$ and $\bar P(s)$ are polynomial matrices with the same eigenstructure, then $\rank(P(s))=\rank(\bar P(s))$ (for instance, they have the same number of homogeneous invariant factors), and   therefore,  by Lemma \ref{theoDeDoMa}, $\deg(P(s))=\deg(\bar P(s))$.

In \cite{AmBaMaRo23} the following row completion   problem was posed:
\begin{problem}\label{problem}
Let $P(s)\in\efe[s]^{m\times n}$ be a polynomial matrix of $\deg(P(s))=d$. Find necessary and sufficient conditions for the existence of a polynomial matrix  $W(s)\in \efe[s]^{z\times n}$ of $\deg(W(s))\leq d$
such that  $\begin{bmatrix}P(s)\\W(s)\end{bmatrix}$  has the eigenstructure, or part 
 of it, prescribed.
\end{problem}
Since the eigenstructure consists of four types of invariants, depending on the invariant(s) prescribed we get different problems. Hence, Problem \ref{problem} leads  actually to 15 different problems, 6 of which  have been  solved in \cite{AmBaMaRo23}. In this paper we solve the remaining cases. 
 By transposition, and interchanging the row and column minimal indices, the  results apply for the column completion problem.

  We  use the solution to Problem \ref{problem} presented in \cite{AmBaMaRo23} when the whole  eigenstructure is prescribed. To state it we need some notation and  definitions.

Given two integers $n$ and $m$, whenever $n>m$ we take  $\sum_{i=n}^{m}=0$. In the same way, if a condition is stated for $n\leq i\leq m$ with $n>m$, we understand that the condition 
 is trivially satisfied.
Let $a_1, \dots,  a_m$ be a  sequence of integers. Whenever we write  $\ba=(a_1, \dots,  a_m)$, we  understand that   $a_1\geq \dots \geq a_m$,
and  we  take $a_i=\infty$ for $i<1$ and $a_i=-\infty$ for $i>m$. If $a_m\geq 0$, the sequence $\ba=(a_1, \dots,  a_m)$ is called a {\em partition}.  

Let   $\ba= (a_1,  \ldots, a_m)$ and $\bb= (b_1, \ldots, b_m)$ be  two sequences of integers. It is said that   $\ba$ is {\em majorized} by $\bb$ (denoted by $\ba \prec \bb$) if $\sum_{i=1}^k a_i \leq \sum_{i=1}^k b_i $ for $1 \leq k \leq m-1$ and $\sum_{i=1}^m a_i =\sum_{i=1}^m b_i$ (this is an extension to sequences of integers of the definition of majorization given for partitions in \cite{HLP88}).

We introduce next the concept of generalized majorization. 
\begin{definition}{\rm \cite[Definition 2]{DoStEJC10}}
Let $\bd = (d_1, \dots, d_m)$, $\ba=(a_1, \dots, a_s)$ and $\bg=(g_1, \dots, g_{m+s})$  be sequences of  integers.
We say that  $\bg$ is majorized by $\bd$ and $\ba$  $(\bg \prec' (\bd,\ba))$ if
\begin{equation}\label{gmaj1}
d_i\geq g_{i+s}, \quad 1\leq i\leq m,
\end{equation}
\begin{equation}\label{gmaj2}
\sum_{i=1}^{h_j}g_i-\sum_{i=1}^{h_j-j}d_i\leq \sum_{i=1}^j a_i, \quad 1\leq j\leq s,
\end{equation}
where $h_j=\min\{i\; : \; d_{i-j+1}<g_i\}$, $1\leq j\leq s$  $(d_{m+1}=-\infty)$,
\begin{equation}\label{gmaj3}
\sum_{i=1}^{m+s}g_i=\sum_{i=1}^md_i+\sum_{i=1}^sa_i.
\end{equation}
\end{definition}
\begin{remark}\label{gmajcp} Observe that if $s=0$, the generalized majorization reduces to   $\bd =\bg$, and if $m=0$ it reduces to  $\bg \prec \ba$.   Also, if $\bg \prec' (\bd,\ba)$ and $\ba\prec \hat \ba$ for some $\hat \ba$, then  $\bg \prec' (\bd,\hat \ba)$.
\end{remark}
 
Given two sequences of integers $\bu = (u_1, \dots, u_{p})$ and $\bb = (b_1, \dots, b_{y})$ the union, $\bu\cup \bb$,  is the decreasingly ordered sequence  of the $p+y$ integers of $\bu$ and $\bb$.
We need the following technical Lemma.
\begin{lemma} {\rm (\cite[Lemma 4.4]{AmBaMaRo23})}\label{lemmacup}
 Let $\bu = (u_1, \dots, u_{p})$ and $\bb = (b_1, \dots, b_{y})$    be sequences of  integers.  
Then
$
\bu\cup \bb\prec'(\bu, \bb).
$
\end{lemma}

The next theorem provides a solution to Problem \ref{problem} when the whole  eigenstructure is prescribed. It is a generalization  to polynomial matrices  of the corresponding result for matrix pencils given in \cite{Do13} (see also \cite{DoSt19}).
\begin{theorem} {\rm(\cite[Theorem 4.2]{AmBaMaRo23}, Prescription of  the  eigenstructure)}
  \label{theoprescr4}
  	Let $P(s)\in\efe[s]^{m\times n}$ be a polynomial matrix, $\deg(P(s))=d$, $\rank (P(s))=r$.
	Let $\phi_1(s,t)\mid\cdots\mid\phi_r(s,t)$ be the
	homogeneous  invariant factors, 
	$\bc=(c_1,  \dots,  c_{n-r})$  the
	column minimal indices, and $\bu=(u_1, \dots,  u_{m-r})$  the row minimal
	indices of $P(s)$, where   $u_1 \geq \dots \geq u_{\eta}  >  u_{\eta +1}= \dots = u_{m-r}=0 $.

        Let $z$ and $x$ be integers such that $0\leq x\leq \min\{z, n-r\}$.  
	Let $\gamma_1(s,t)\mid\cdots\mid\gamma_{r+x}(s,t)$ be monic homogeneous polynomials, and
	$\bd=(d_1, \dots, d_{n-r-x})$
	 and 
	$\bv=(v_1, \dots, v_{m+z-r-x})$ two partitions, where 
	$v_1 \geq  \dots  \geq  v_{\bar \eta}  > v_{\bar \eta +1} = \dots =v_{m+z-r-x} = 0$.
	There exists a polynomial matrix  $W(s)\in \efe[s]^{z\times n}$ such that $\deg(W(s))\leq d$, $\rank \left(\begin{bmatrix}P(s)\\W(s)\end{bmatrix}\right) =r+x$, and 
	$\begin{bmatrix}P(s)\\W(s)\end{bmatrix}$ has $\gamma_1(s,t)\mid\cdots\mid\gamma_{r+x}(s,t)$ as homogeneous invariant factors,
	$d_1, \dots, d_{n-r-x}$ as column minimal indices
	and 
	$v_1, \dots, v_{m+z-r-x}$ as row minimal indices  
	if and only if 
	\begin{equation}\label{eqinterfipolhom}
		\gamma_i(s,t)\mid \phi_i(s,t)\mid \gamma_{i+z}(s,t),\quad 1\leq i \leq r,
	\end{equation}
	\begin{equation}\label{eqetapol}\bar \eta \geq  \eta,\end{equation}
	\begin{equation}\label{eqcmimajpol}  \bc \prec'  (\bd , \ba),\end{equation}
	\begin{equation}\label{eqrmimajpol}\bv \prec'  (\bu , \bb),\end{equation}
        \begin{equation}\label{eqdegsumpol}
 \begin{aligned}
		\sum_{i=1}^{r+x}\deg(\lcm(\phi_{i-x},\gamma_i))
		\leq \sum_{i=1}^{m+z-r-x}v_i-\sum_{i=1}^{m-r}u_i+\sum_{i=1}^{r+x}\deg(\gamma_i),
  \\
\mbox{ with equality when $x=0$,}
\end{aligned}
\end{equation}
	where $\ba = (a_1, \dots, a_x )$ and $\bb = (b_1, \dots,  b_{z-x} )$ are
	\begin{equation}\label{eqdefa}
		\begin{array}{rl}
			a_1=&
			\sum_{i=1}^{m+z-r-x}v_i-\sum_{i=1}^{m-r} u_i\\&+\sum_{i=1}^{r+x}\deg( \gamma_i)-
			\sum_{i=1}^{  r+x-1}\deg(\lcm(  \phi_{i-x+1},  \gamma_i))-d,\\
			a_j=&
			\sum_{i=1}^{  r+x-j+1}\deg(\lcm(  \phi_{i-x+j-1},  \gamma_i))
			-
			\sum_{i=1}^{  r+x-j}\deg(\lcm(  \phi_{i-x+j},  \gamma_i))
			-d,\\& \hfill 2\leq j \leq x,
		\end{array}
	\end{equation}
	\begin{equation}\label{eqdefb}
		\begin{array}{rl}
			b_1=&	   \sum_{i=1}^{m+z-r-x}  v_i-\sum_{i=1}^{m-r}  u_i\\&+\sum_{i=1}^{  r+x}\deg(  \gamma_i)-
			\sum_{i=1}^{  r+x}\deg(\lcm(  \phi_{i-x-1},  \gamma_i)),
			\\
			b_j=&
			\sum_{i=1}^{  r+x}\deg(\lcm(  \phi_{i-x-j+1},  \gamma_i)-
			\sum_{i=1}^{  r+x}\deg(\lcm(  \phi_{i-x-j},  \gamma_i)),\\&\hfill  2\leq j \leq z-x.
		\end{array}
	\end{equation}

\end{theorem}
\begin{remark}\label{remposttheo}
  \begin{enumerate}
  \item\label{remitdecreasing}  
 
  From \cite[Lemmas 1 and 2]{DoStEJC10} (see also \cite[Remark 2.7]{AmBaMaRo23}) we derive  
  $  a_1\geq \dots\geq   a_x$ and $  b_1\geq \dots\geq   b_{z-x}\geq 0$.
  
\item\label{remituvtocd}
It was seen in   \cite[Remark 4.3]{AmBaMaRo23}
that if 
(\ref{eqinterfipolhom})-(\ref{eqdegsumpol}) hold, then
\begin{equation}\label{eqdegsumpolequiv}
\begin{aligned}
  \sum_{i=1}^{r+x}\deg(\lcm(\phi_{i-x},\gamma_i))
  \leq \sum_{i=1}^{n-r}c_i-\sum_{i=1}^{n-r-x}d_i+\sum_{i=1}^{r}\deg(\phi_i)+xd,\\
  \mbox{with equality when $x=z$,}
  \end{aligned}
\end{equation}
and (\ref{eqcmimajpol}) and  (\ref{eqrmimajpol}) hold for
   $\ba=(a_1, \dots, a_x)$ and $\bb=(b_1, \dots, b_{z-x})$ defined as
\begin{equation}\label{eqdefabis}
\begin{array}{rl}
a_1=&
  \sum_{i=1}^{n-r}c_i-\sum_{i=1}^{n-r-x}d_i\\&+\sum_{i=1}^{r}\deg(\phi_i)
  -
\sum_{i=1}^{  r+x-1}\deg(\lcm(  \phi_{i-x+1},  \gamma_i))+(x-1)d,\\
a_j=&
\sum_{i=1}^{  r+x-j+1}\deg(\lcm(  \phi_{i-x+j-1},  \gamma_i))\\
&-
\sum_{i=1}^{  r+x-j}\deg(\lcm(  \phi_{i-x+j},  \gamma_i))
-d,\quad   2\leq j \leq x,
\end{array}
\end{equation}
\begin{equation}\label{eqdefbbis}
\begin{array}{rl}
 b_1=&
\sum_{i=1}^{n-r}c_i-\sum_{i=1}^{n-r-x}d_i+\sum_{i=1}^{r}\deg(\phi_i)\\&-
\sum_{i=1}^{  r+x}\deg(\lcm(  \phi_{i-x-1},  \gamma_i))+xd,
\\
  b_j=&
\sum_{i=1}^{  r+x}\deg(\lcm(  \phi_{i-x-j+1},  \gamma_i)-
\sum_{i=1}^{  r+x}\deg(\lcm(  \phi_{i-x-j},  \gamma_i)),\\& \hfill 2\leq j \leq z-x,
\end{array}
\end{equation}

Conversely, (\ref{eqinterfipolhom})-(\ref{eqrmimajpol}) and (\ref{eqdegsumpolequiv}) 
with $\ba$ and $\bb$   defined as in 
(\ref{eqdefabis}) and (\ref{eqdefbbis}), respectively, imply (\ref{eqinterfipolhom})-(\ref{eqdegsumpol}),
 with $\ba$ and $\bb$   defined as in 
(\ref{eqdefa}) and (\ref{eqdefb}).

  \end{enumerate}
 \end{remark}

 \section{Row (column) completion for polynomial matrices: remaining cases}
\label{seccompl}

The aim of this section is to give a solution to  Problem \ref{problem}  for the cases  not 
solved in \cite{AmBaMaRo23}. We present the solution in different subsections depending on the invariant(s) prescribed.
The invariants that are not prescribed have a certain degree of freedom, but they must satisfy some restrictions determined by Theorem \ref{theoprescr4}
  (see Examples \ref{exalgclosed} and \ref{exiedcmi} below).

\subsection{Prescription of infinite and singular structures} 
We start  by   prescribing the infinite elementary divisors and the column and row minimal indices.  
\begin{theorem} {\rm (Prescription of infinite and singular structures)}
  \label{theoprescrdeirmicmi}
  Let $\efe$ be an algebraically closed field.
Let $P(s)\in\efe[s]^{m\times n}$,  $\deg(P(s))=d$, $\rank (P(s))=r$.
  Let $\alpha_1(s)\mid\cdots\mid\alpha_r(s)$ be the invariant factors,
$e_1\leq\cdots\leq e_r$
the partial multiplicities of $\infty$, 
$\bc=(c_1,  \dots,  c_{n-r})$  the
column minimal indices, and $ \bu=(u_1, \dots,  u_{m-r})$  the row minimal
indices of $P(s)$, where 
$u_1 \geq  \dots  \geq  u_{\eta}  > u_{\eta +1} = \dots =u_{m-r} = 0$.

Let $z$ and $x$ be integers such that $0\leq x\leq \min\{z, n-r\}$ and 
let  $f_1\leq\dots \leq f_{r+x} $,
 $d_1 \geq  \dots  \geq  d_{n-r-x}\geq 0$ and
 $v_1 \geq  \dots  \geq  v_{\bar \eta} > v_{\bar \eta+1}=\dots= v_{m+z-r-x}=0$ be non negative integers.
Let
\begin{equation}\label{eqbigAdeicmirmi}
A=\sum_{i=1}^{r+x}f_i-\sum_{i=1}^{r}e_i+\sum_{i=1}^{n-r-x}d_i-\sum_{i=1}^{n-r}c_i+\sum_{i=1}^{m+z-r-x}v_i-\sum_{i=1}^{m-r}u_i-xd.
\end{equation}
There exists $W(s)\in \efe[s]^{z\times n}$ such that $\deg(W(s))\leq d$, $\rank \left(\begin{bmatrix}P(s)\\W(s)\end{bmatrix}\right)=r+x$
and $\begin{bmatrix}P(s)\\W(s)\end{bmatrix}$  has 
 $f_1\leq\dots \leq f_{r+x} $ as  partial multiplicities of $\infty$, 
 $\bd=(d_1,\dots, d_{n-r-x})$ as  column minimal indices and
$\bv=(v_1,\dots, v_{m+z-r-x})$ as row minimal indices
if and only if 
(\ref{eqetapol}),
\begin{equation}\label{eqinteriedpol}
 f_i\leq e_i \leq f_{i+z},\quad 1\leq i \leq r,
  \end{equation}
\begin{equation}\label{eqAleq}
A\leq \sum_{i=r+x-z+1}^{r}\deg(\alpha_i),
\end{equation}
\begin{equation}\label{eqvusA}
 \sum_{i=1}^{m+z-r-x}v_i-\sum_{i=1}^
 {m-r}u_i\geq \max\{0, A\}+\sum_{i=1}^{r+x}\max\{e_{i-x}, f_i\}-\sum_{i=1}^{r+x}f_i,
\end{equation}
\begin{equation}\label{eqcmimajpolhatsA}  \bc \prec'  (\bd , \hat \ba),\end{equation}
\begin{equation}\label{eqrmimajpolhatsA} \bv \prec'  (\bu , \hat \bb),\end{equation}
where $\hat \ba = (\hat a_1, \dots, \hat a_x )$ and $\hat \bb = (\hat b_1, \dots,  \hat b_{z-x})$ 
are 
\begin{equation}\label{eqadeicmirmi}\begin{array}{rl}\hat a_1=&
  \sum_{i=1}^{m+z-r-x}v_i-\sum_{i=1}^{m-r}u_i\\&+\sum_{i=1}^{r+x}f_i-\sum_{i=1}^{r+x-1}\max\{e_{i-x+1}, f_i\}-A-d
  ,\\ \hat a_j=&
\sum_{i=1}^{r+x-j+1}\max\{e_{i-x+j-1}, f_i\}\\&-\sum_{i=1}^{r+x-j}\max\{e_{i-x+j}, f_i\}
  -d,\quad 2\leq j \leq x,
\end{array}
\end{equation}
\begin{equation}\label{eqbdeicmirmi}
\begin{array}{rl}\hat b_1=& \sum_{i=1}^{m+z-r-x}v_i-\sum_{i=1}^{m-r}u_i+\min\{0,\deg(\alpha_{r})-A \}\\&
  +\sum_{i=1}^{  r+x}f_i-\sum_{i=1}^{  r+x}\max\{e_{i-x-1}, f_i\},
  \\ \hat b_j=&\min\{0,\sum_{i=r-j+1}^{r}\deg(\alpha_{i})-A\}-\min\{0,\sum_{i=r-j+2}^{r}\deg(\alpha_{i})-A \}\\&
 +\sum_{i=1}^{  r+x}\max\{e_{i-x-j+1}, f_i\}-\sum_{i=1}^{  r+x}\max\{e_{i-x-j}, f_i\}
  ,\quad 2\leq j \leq z-x.
\end{array}
\end{equation}
\end{theorem}
\begin{remark}\label{remacdecr}\ 
  \begin{enumerate}
    \item\label{remitac}
The necessity of the conditions,  and the sufficiency when $A\leq0$ hold for arbitrary fields.
  \item\label{remitemhatadecr} 
   If (\ref{eqvusA}) holds, then,  as in   Remark \ref{remposttheo}.\ref{remitdecreasing}, we have $\hat a_1\geq \dots \geq \hat a_x$.
     \item\label{remitemhatbdecr}   
        Let us see that  if (\ref{eqAleq}) and (\ref{eqvusA}) hold then  $\hat b_1\geq \dots \geq \hat b_{z-x}$.
        By (\ref{eqAleq}), $z-x\in \{k\geq 0\; : \;A\leq \sum_{i=r-k+1}^{r}\deg(\alpha_{i}) \}$. Let $g=\min\{k\geq 0\; : \;A\leq \sum_{i=r-k+1}^{r}\deg(\alpha_{i}) \}$.
        Then $0\leq g \leq z-x$. Observe that  $g=0$ if and only if $A\leq 0$. If $g\geq 1$ then
        \begin{equation}\label{eqtarte}
          \sum_{i=r-g+2}^{r}\deg(\alpha_{i}) <A\leq \sum_{i=r-g+1}^{r}\deg(\alpha_{i}).\end{equation}
        If $g=0$ or $g=1$ then $A\leq \deg(\alpha_r)\leq  \sum_{i=r-j+1}^{r}\deg(\alpha_{i}) $ for $1\leq j \leq z-x$; hence
        $$\begin{array}{rl}\hat b_1=& \sum_{i=1}^{m+z-r-x}v_i-\sum_{i=1}^{m-r}u_i +\sum_{i=1}^{  r+x}f_i-\sum_{i=1}^{  r+x}\max\{e_{i-x-1}, f_i\},       
  \\ \hat b_j=& \sum_{i=1}^{  r+x}\max\{e_{i-x-j+1}, f_i\}-\sum_{i=1}^{  r+x}\max\{e_{i-x-j}, f_i\},\\&\hfill 2\leq j \leq z-x,
\end{array}
        $$
        and, as in  Remark  \ref{remposttheo}.\ref{remitdecreasing}, we have $\hat b_1\geq \dots \geq \hat b_{z-x}\geq 0.$
       
        If $g\geq 2$, then from (\ref{eqtarte}) we get
        $$
        \begin{array}{l}
        \sum_{i=r-j+2}^{r}\deg(\alpha_{i})\leq\sum_{i=r-j+1}^{r}\deg(\alpha_{i})\leq \sum_{i=r-g+2}^{r}\deg(\alpha_{i}) <A, \\\hfill 1\leq j \leq g-1,\\
        A\leq \sum_{i=r-g+1}^{r}\deg(\alpha_{i})\leq \sum_{i=r-j+2}^{r}\deg(\alpha_{i})\leq \sum_{i=r-j+1}^{r}\deg(\alpha_{i}),\\\hfill g+1\leq j \leq z-x.
        \end{array}
        $$
        Thus, 
$$\begin{array}{rl}\hat b_1=&  \sum_{i=1}^{m+z-r-x}v_i-\sum_{i=1}^{m-r}u_i-A\\&+\deg(\alpha_r)+\sum_{i=1}^{  r+x}f_i-\sum_{i=1}^{  r+x}\max\{e_{i-x-1}, f_i\},
  \\ \hat b_j=& \deg(\alpha_{r-j+1})+\sum_{i=1}^{  r+x}\max\{e_{i-x-j+1}, f_i\}\\&-\sum_{i=1}^{  r+x}\max\{e_{i-x-j}, f_i\},\quad 2\leq j \leq g-1,
  \\\hat b_{g}=&A-\sum_{i=r-g+2}^{r}\deg(\alpha_{i})+\sum_{i=1}^{  r+x}\max\{e_{i-x-g+1}, f_i\}\\&-\sum_{i=1}^{  r+x}\max\{e_{i-x-g}, f_i\},
  \\\hat b_{j}=&\sum_{i=1}^{  r+x}\max\{e_{i-x-j+1}, f_i\}-\sum_{i=1}^{  r+x}\max\{e_{i-x-j}, f_i\}, \\& \hfill g+1\leq j\leq z-x.
\end{array}
$$
        From (\ref{eqvusA}) and (\ref{eqtarte}), we get       
        $\sum_{i=1}^{m+z-r-x}v_i-\sum_{i=1}^{m-r}u_i-A+\sum_{i=1}^{  r+x}f_i\geq\sum_{i=1}^{  r+x}\max\{e_{i-x}, f_i\}$ and
        $0<A-\sum_{i=r-g+2}^{r}\deg(\alpha_{i})\leq  \deg(\alpha_{r-g+1})$, respectively.
       Thus, as in 
       Remark \ref{remposttheo}.\ref{remitdecreasing},   we have  $\hat b_1\geq \dots \geq \hat b_{z-x}\geq 0$.
  \end{enumerate}
\end{remark}
{\bf Proof of Theorem \ref{theoprescrdeirmicmi}.}
Let $\phi_1(s,t)\mid \cdots \mid \phi_r(s,t)$ be the homogeneous invariant factors of $P(s)$.

Assume that there exists
 $W(s)\in \efe[s]^{z\times n}$,
$\deg(W(s))\leq d$, $\rank \left(\begin{bmatrix}P(s)\\W(s)\end{bmatrix}\right) =r+x$
and $\begin{bmatrix}P(s)\\W(s)\end{bmatrix}$  has 
 $f_1\leq\dots \leq f_{r+x} $ as  partial multiplicities of $\infty$, 
 $\bd=(d_1,\dots, d_{n-r-x})$ as  column minimal indices and
$\bv=(v_1,\dots, v_{m+z-r-x})$ as row minimal indices.
Let $\beta_1(s)\mid\cdots\mid\beta_{r+x}(s)$ and $\gamma_1(s, t)\mid\cdots\mid\gamma_{r+x}(s,t)$ be its invariant factors and  homogeneous invariant factors, respectively.
Then, by Lemma \ref{theoDeDoMa} we obtain
\begin{equation}\label{eqbetaIST}
A=\sum_{i=1}^{r}\deg(\alpha_i)-\sum_{i=1}^{r+x}\deg(\beta_i),
\end{equation}
and
by Theorem \ref{theoprescr4}, (\ref{eqinterfipolhom})-(\ref{eqdegsumpol}) hold, where $\ba$ and $\bb$ are defined in (\ref{eqdefa}) and (\ref{eqdefb}), respectively.
Condition (\ref{eqinterfipolhom}) is equivalent to  (\ref{eqinteriedpol}) and
\begin{equation}\label{eqinterfifpol}
\beta_i(s) \mid \alpha_i(s) \mid \beta_{i+z}(s), \quad 1\leq i \leq r.  
\end{equation}
Taking into account that
$$
\begin{array}{rl}
&\sum_{i=1}^{  r+x}\deg(\lcm(  \alpha_{i-x},  \beta_i))-\sum_{i=1}^{r+x}\deg( \beta_i)
\\\geq&
\max \{\sum_{i=1}^{r}\deg( \alpha_{i})-\sum_{i=1}^{r+x}\deg( \beta_i), 0\}=\max\{A,0\},\\
\end{array}
$$
from (\ref{eqdegsumpol}) we obtain (\ref{eqvusA}).
From (\ref{eqinterfifpol}) we have $\sum_{i=1}^{r+x}\deg(\beta_i)=\sum_{i=1}^{z}\deg(\beta_i)+\sum_{i=1}^{r+x-z}\deg(\beta_{i+z})\geq \sum_{i=1}^{r+x-z}\deg(\alpha_{i});$
hence we obtain (\ref{eqAleq}). 
We also have
$$\begin{array}{rl}\sum_{i=1}^{x}\hat a_i= &\sum_{i=1}^{m+z-r-x}v_i-\sum_{i=1}^{m-r}u_i+\sum_{i=1}^{r+x}f_i-\sum_{i=1}^{r}e_i-A-xd\\=&\sum_{i=1}^{n-r}c_i-
\sum_{i=1}^{n-r-x}d_i=\sum_{i=1}^{x}a_i,\end{array}$$
$$
\begin{array}{rl}
  \sum_{i=j}^{x}a_i=&\sum_{i=1}^{r+x-j+1}\deg(\lcm(  \phi_{i-x+j-1},  \gamma_i))-\sum_{i=1}^{r}\deg(\phi_i)\\&-(x-j+1)d\\
  =&\sum_{i=1}^{r+x-j+1}\max\{e_{i-x+j-1}, f_i\}-\sum_{i=1}^{r}e_i-(x-j+1)d\\&+\sum_{i=1}^{r+x-j+1}\deg(\lcm(  \alpha_{i-x+j-1},  \beta_i))-\sum_{i=1}^{r}\deg(\alpha_i)\\=&
  \sum_{i=j}^x \hat a_i+\sum_{i=1}^{r+x-j+1}\deg(\lcm(  \alpha_{i-x+j-1},  \beta_i))-\sum_{i=1}^{r}\deg(\alpha_i)\\\geq&\sum_{i=j}^x \hat a_i, \quad 2\leq j\leq x,
\end{array}
$$
and
$$\sum_{i=1}^{j}a_i=\sum_{i=1}^{x}a_i-\sum_{i=j+1}^{x}a_i\leq \sum_{i=1}^{x}\hat a_i-\sum_{i=j+1}^{x}\hat a_i, \quad 1\leq j\leq x-1.
$$
 Therefore, $\ba\prec \hat \ba$, and  from (\ref{eqcmimajpol}) 
 and Remark \ref{gmajcp} 
 we derive  (\ref{eqcmimajpolhatsA}). 
From (\ref{eqAleq}), (\ref{eqinteriedpol}) and (\ref{eqinterfipolhom}) we get
$\sum_{i=1}^{z-x}\hat b_i=\sum_{i=1}^{m+z-r-x}v_i-\sum_{i=1}^{m-r}u_i=
\sum_{i=1}^{z-x} b_i.$
 Moreover,
$$
\begin{array}{rl}
  \sum_{i=j}^{z-x}\hat b_i=&-\min\{0, \sum_{i=r-j+2}^{r}\deg(\alpha_{i})-A \}
 \\& +\sum_{i=1}^{r+x}\max\{e_{i-x-j+1}, f_i\}-\sum_{i=1}^{r+x}f_i
  \\=&
  \max\{0, A-\sum_{i=r-j+2}^{r}\deg(\alpha_{i})\}
  \\&+
  \sum_{i=1}^{r+x}\max\{e_{i-x-j+1}, f_i\}-\sum_{i=1}^{r+x}f_i\\
  =&\max\{0,\sum_{i=1}^{r-j+1}\deg(\alpha_{i})-\sum_{i=1}^{r+x}\deg(\beta_{i}) \}
  \\&+\sum_{i=1}^{r+x}\max\{e_{i-x-j+1}, f_i\}-\sum_{i=1}^{r+x}f_i, \quad 2\leq j\leq z-x,
\end{array}
$$
$$
\begin{array}{rl}
  \sum_{i=j}^{z-x}b_i=&
  \sum_{i=1}^{r+x}\deg(\lcm(\alpha_{i-x-j+1},\beta_i))-\sum_{i=1}^{r+x}\deg(\beta_i)
  \\&+ \sum_{i=1}^{r+x}\max\{e_{i-x-j+1},f_i\}-
  \sum_{i=1}^{r+x}f_i\\
  \geq &\max\{\sum_{i=1}^{r-j+1}\deg(\alpha_{i})-\sum_{i=1}^{r+x}\deg(\beta_{i}), 0 \}
  \\&+\sum_{i=1}^{r+x}\max\{e_{i-x-j+1},f_i\}-
  \sum_{i=1}^{r+x}f_i
  \\=&
 \sum_{i=j}^{z-x}\hat b_i,
\quad 2\leq j\leq z-x, 
  \end{array}
$$
and
$$\sum_{i=1}^{j}b_i=\sum_{i=1}^{z-x}b_i-\sum_{i=j+1}^{z-x}b_i\leq \sum_{i=1}^{z-x}\hat b_i-\sum_{i=j+1}^{z-x}\hat b_i, \quad 1\leq j\leq z-x-1.
$$
Thus, $\bb\prec \hat \bb$, and   from (\ref{eqrmimajpol}) 
 and Remark \ref{gmajcp}
we derive  (\ref{eqrmimajpolhatsA}).

Conversely, 
assume that (\ref{eqetapol}) and
(\ref{eqinteriedpol})-(\ref{eqrmimajpolhatsA}) hold.  
We will prove that there exist monic polynomials $\beta_1(s)\mid\cdots\mid\beta_{r+x}(s)$ such that, if
$
\gamma_i(s,t)=t^{f_i}t^{\deg(\beta_i)}\beta_i(\frac{s}{t}), 
$
$1\leq i \leq r+x$, then (\ref{eqinterfipolhom}) and 
(\ref{eqcmimajpol})-(\ref{eqdegsumpol}) are satisfied.
Once the existence of these polynomials is proven, 
 by Theorem \ref{theoprescr4}, there exists  $W(s)\in \efe[s]^{z\times n}$
 such that $\deg(W(s))\leq d$, $\rank \left(\begin{bmatrix}P(s)\\W(s)\end{bmatrix}\right) =r+x$, and 
	$\begin{bmatrix}P(s)\\W(s)\end{bmatrix}$ has $\gamma_1(s,t)\mid\cdots\mid\gamma_{r+x}(s,t)$ as homogeneous invariant factors,
	$\bd=(d_1, \dots, d_{n-r-x})$ as column minimal indices
	and 
	$\bv=(v_1, \dots, v_{m+z-r-x})$ as row minimal indices. By the definition of $\gamma_1(s,t),\dots,\gamma_{r+x}(s,t)$, the exponents  $f_1,\dots,f_{r+x}$ are the partial multiplicities of $\infty$ in
$\begin{bmatrix}P(s)\\W(s)\end{bmatrix}$.

Therefore, we only need to prove the existence of   polynomials $\beta_1(s)\mid\cdots\mid\beta_{r+x}(s)$   satisfying 
the conditions mentioned above.

Observe that if $x=0$ then, from (\ref{eqcmimajpolhatsA}), $\bc =\bd$  (see Remark \ref{gmajcp});  hence $A=\sum_{i=1}^{r}f_i-\sum_{i=1}^{r}e_i+\sum_{i=1}^{m+z-r}v_i-\sum_{i=1}^{m-r}u_i$ and from
(\ref{eqvusA}) we obtain $A\geq 0$.
If $x=z$ then, from (\ref{eqAleq}), we obtain $A\leq 0$. Hence, we distinguish two cases.

\begin{itemize}
  \item  
Let $A>0$.  We saw  that if $g=\min\{k\geq 0\; : \;A\leq \sum_{i=r-k+1}^{r}\deg(\alpha_{i}) \}$ and    
 $A>0$, then $1\leq g\leq z-x$ (see Remark \ref{remacdecr}). Therefore  (\ref{eqtarte}) holds, i.e.,    
$0<A-\sum_{i=r-g+2}^{r}\deg(\alpha_{i})\leq \deg(\alpha_{r-g+1})
.$

Let $h=\min\{k\; : \;A-\sum_{i=r-g+2}^{r}\deg(\alpha_{i})\leq \deg(\alpha_{k-g+1}) \}$. Obviously, $h\leq r$.
For $j\leq g-1$, we have 
$\deg(\alpha_{j-g+1})=0<A-\sum_{i=r-g+2}^{r}\deg(\alpha_{i})$; hence $h\geq g$
and
    \begin{equation}\label{eqdeg}\deg(\alpha_{h-g})< A-\sum_{i=r-g+2}^{r}\deg(\alpha_{i})\leq \deg(\alpha_{h-g+1}).\end{equation}
    Taking $$w=
    \deg(\alpha_{h-g})+\deg(\alpha_{h-g+1})+\sum_{i=r-g+2}^{r}\deg(\alpha_{i})-A,$$
     we get
    $\deg(\alpha_{h-g})\leq w<\deg(\alpha_{h-g+1})$.
    Since $\efe$ is an algebraically closed field, there exists a monic polynomial
    $\tau(s)$  such that
$$\alpha_{h-g}(s)\mid \tau(s)\mid \alpha_{h-g+1}(s), \quad \deg(\tau)=w.$$
Define 
$$
\begin{array}{rll}
  \beta_i(s)=&\alpha_{i-x-g}(s),&1\leq i\leq h+x-1,
  \\
  \beta_{h+x}(s)=&\tau(s),\\
  \beta_i(s)=&\alpha_{i-x-g+1}(s),&h+x+1\leq i\leq r+x.
\end{array}
$$
We have
$\beta_1(s)\mid \dots \mid\beta_{r+x}(s)$, and 
$$\begin{array}{rl}\sum_{i=1}^{r+x}\deg(\beta_{i})=&\sum_{i=1}^{h-g+1}\deg(\alpha_{i})+\sum_{i=r-g+2}^{r}\deg(\alpha_{i})-A\\&+\sum_{i=h-g+2}^{r-g+1}\deg(\alpha_{i}),
\end{array}$$
 i.e., (\ref{eqbetaIST}) holds.

 Moreover,
$\beta_i(s)\mid \alpha_{i-x-g+1}(s)\mid\alpha_{i}(s)\mid\beta_{i+x+g}(s)\mid\beta_{i+z}(s)$ for $1\leq i \leq r$, 
then (\ref{eqinterfifpol}) holds. From (\ref{eqinteriedpol}) and (\ref{eqinterfifpol}) we obtain
(\ref{eqinterfipolhom}).
Since $\beta_i(s)\mid \alpha_{i-x-g+1}(s)\mid \alpha_{i-x+j}(s)$, $1\leq i \leq r+x-j$, $0\leq j \leq x$, from (\ref{eqbetaIST})  it follows that
$$\sum_{i=1}^{r+x-j}\deg(\lcm(\alpha_{i-x+j}, \beta_i))=\sum_{i=1}^{r}\deg(\alpha_{i})=A+\sum_{i=1}^{r+x}\deg(\beta_i), \; 0\leq j \leq x;$$
hence, for $0\leq j \leq x$,
$$\begin{array}{ll}
  &\sum_{i=1}^{r+x-j}\deg(\lcm(\phi_{i-x+j}, \gamma_i))-\sum_{i=1}^{r+x}\deg(\gamma_i)\\=&A+\sum_{i=1}^{r+x-j}\max\{e_{i-x+j}, f_i\}-\sum_{i=1}^{r+x}f_i.\end{array}$$
If $x=0$, recalling  that  $A=\sum_{i=1}^{r}f_i-\sum_{i=1}^{r}e_i+\sum_{i=1}^{m+z-r}v_i-\sum_{i=1}^{m-r}u_i$, from (\ref{eqinteriedpol}) we obtain
  $\sum_{i=1}^{r}\deg(\lcm(\phi_{i}, \gamma_i))= \sum_{i=1}^{m+z-r}v_i-\sum_{i=1}^{m-r}u_i+ \sum_{i=1}^{r}\deg(\gamma_i)$. 
Thus, (\ref{eqvusA}) implies   (\ref{eqdegsumpol}), and if $\ba$ is defined as in (\ref{eqdefa}), then
$$\begin{array}{rl}
  \sum_{i=1}^{j}a_i=&
\sum_{i=1}^{m+z-r-x}v_i-\sum_{i=1}^{m-r} u_i\\&+\sum_{i=1}^{r+x}f_i-
\sum_{i=1}^{  r+x-j}\max\{e_{i-x+j},  f_i\}-A-jd\\
=&\sum_{i=1}^{j}\hat a_i, \quad 1\leq j\leq x,
\end{array}
$$
 therefore $\ba=\hat \ba$ and (\ref{eqcmimajpolhatsA}) is equivalent to (\ref{eqcmimajpol}).

Let $\bb$ be as in (\ref{eqdefb}).
Let $j\in\{1,\dots, g-1\}$. Then $\sum_{i=r-j+1}^{r}\deg(\alpha_{i})
<A$  and $\beta_i(s)\mid\alpha_{i-x-g+1}(s)\mid\alpha_{i-x-j}(s)$, $1\leq i\leq r+x$; hence
$$\sum_{i=1}^{  r+x}\deg(\lcm(  \alpha_{i-x-j},  \beta_i))=\sum_{i=1}^{r-j}\deg(\alpha_i),$$ and, taking into account (\ref{eqbetaIST}), 
we obtain
$$\begin{array}{rl}
\sum_{i=1}^{j}b_i
=&
\sum_{i=1}^{m+z-r-x}  v_i-\sum_{i=1}^{m-r}  u_i-A+\sum_{i=r-j+1}^{  r}\deg( \alpha_i)\\&+\sum_{i=1}^{  r+x}f_i-\sum_{i=1}^{  r+x}\max\{e_{i-x-j}, f_i\}=\sum_{i=1}^{j}\hat b_i.
\end{array}
$$
Let $j\in\{g,\dots, z-x\}$. Then $A\leq \sum_{i=r-j+1}^{r}\deg(\alpha_{i})$ and $$\alpha_{i-x-j}(s)\mid\alpha_{i-x-g}(s)\mid\beta_{i}(s),\quad 1\leq i\leq r+x;$$ hence
$$\sum_{i=1}^{  r+x}\deg(\lcm(  \alpha_{i-x-j},  \beta_i))=\sum_{i=1}^{  r+x}\deg(\beta_i),$$ and
$$
\begin{array}{rl}
\sum_{i=1}^{j}b_i
=&
\sum_{i=1}^{m+z-r-x}  v_i-\sum_{i=1}^{m-r}  u_i
+\sum_{i=1}^{  r+x}f_i\\&-\sum_{i=1}^{  r+x}\max\{e_{i-x-j}, f_i\}=
\sum_{i=1}^{j}\hat b_i.
\end{array}
$$
Therefore $\bb=\hat \bb$ and (\ref{eqrmimajpolhatsA}) is equivalent to (\ref{eqrmimajpol}).
\item  
  If  $A\leq 0$, 
  let $\tau(s)$ be a monic polynomial such that $\deg(\tau)=-A$.
  Define
$$
\begin{array}{rll}
  \beta_i(s)=&\alpha_{i-x}(s),&1\leq i\leq r+x-1,
  \\
  \beta_{r+x}(s)=&\alpha_{r}(s)\tau(s).\\
\end{array}
$$
Then $\beta_1(s)\mid \dots \mid \beta_{r+x}(s)$ and  (\ref{eqbetaIST}) holds.

 If $x=0$, recall that $A\geq 0$. Then, $A=0$, $\tau(s)=1$, 
$\beta_i(s)=\alpha_i(s)$, $1\leq i \leq r$, and (\ref{eqinterfifpol}) holds.

If $x\geq 1$ we have
$$\beta_i(s)=\alpha_{i-x}(s)\mid \alpha_i(s), \quad 1\leq i \leq r,$$
$$\alpha_{i}(s)\mid\alpha_{i+z-x}(s)=\beta_{i+z}(s), \quad 1\leq i \leq r+x-z-1,$$
$$\alpha_{r+x-z}(s)\mid\alpha_{r}(s)\mid \alpha_{r}(s)\tau(s)=\beta_{r+x}(s),$$
$$\alpha_{i}(s)\mid\beta_{i+z}(s)=0, \quad r+x-z< i \leq r,$$
 therefore (\ref{eqinterfifpol}) also holds.
From (\ref{eqinteriedpol}) and (\ref{eqinterfifpol}) we obtain
(\ref{eqinterfipolhom}).

 By definition, 
 $\alpha_{i-x}(s)\mid \beta_{i}(s)$, $1\leq i \leq r+ x$, then $$\sum_{i=1}^{r+x}\deg(\lcm(\phi_{i-x}, \gamma_i))-\sum_{i=1}^{r+x}\deg(\gamma_i)= \sum_{i=1}^{r+x}\max\{e_{i-x}, f_i\}-\sum_{i=1}^{r+x}f_i.$$
 If $x=0$ then $0=A=\sum_{i=1}^{r}f_i-\sum_{i=1}^{r}e_i+\sum_{i=1}^{m+z-r}v_i-\sum_{i=1}^{m-r}u_i$ and from (\ref{eqinteriedpol}) we obtain
  $\sum_{i=1}^{r}\deg(\lcm(\phi_{i}, \gamma_i))= \sum_{i=1}^{m+z-r}v_i-\sum_{i=1}^{m-r}u_i+ \sum_{i=1}^{r}\deg(\gamma_i)$.
Thus, from (\ref{eqvusA}) we obtain (\ref{eqdegsumpol}).

Let
$\ba$ and  $\bb$ be  defined as in (\ref{eqdefa}) and (\ref{eqdefb}), respectively.

Let
$j\in\{1,\dots, x \}$. Then $\beta_i(s)=\alpha_{i-x}(s)\mid\alpha_{i-x+j}(s)$, $1\leq i\leq r+x-j$; hence we obtain
$\sum_{i=1}^{  r+x-j}\deg(\lcm(\alpha_{i-x+j}, \beta_i))=\sum_{i=1}^{  r}\deg(\alpha_{i})$ and, taking into account (\ref{eqbetaIST}),
$$\begin{array}{rl}
  \sum_{i=1}^{j}a_i=&
\sum_{i=1}^{m+z-r-x}v_i-\sum_{i=1}^{m-r} u_i+\sum_{i=1}^{r+x}\deg( \beta_i)-\sum_{i=1}^{  r}\deg(\alpha_{i})\\&+\sum_{i=1}^{r+x}f_i-
\sum_{i=1}^{  r+x-j}\max\{e_{i-x+j},  f_i\}-jd\\
=&
\sum_{i=1}^{m+z-r-x}v_i-\sum_{i=1}^{m-r} u_i-A\\&+\sum_{i=1}^{r+x}f_i-
\sum_{i=1}^{  r+x-j}\max\{e_{i-x+j},  f_i\}-jd
=\sum_{i=1}^{j}\hat a_i.
\end{array}
$$
Therefore $\ba=\hat \ba$ and (\ref{eqcmimajpolhatsA}) is equivalent to (\ref{eqcmimajpol}).

Let $j\in\{1,\dots, z-x\}$. Then $\alpha_{i-x-j}(s)\mid\alpha_{i-x}(s)\mid\beta_{i}(s)$, $1\leq i\leq r+x$; hence we obtain
$\sum_{i=1}^{  r+x}\deg(\lcm(  \alpha_{i-x-j},  \beta_i))=\sum_{i=1}^{  r+x}\deg(\beta_i)$ and
$$
\begin{array}{rl}
\sum_{i=1}^{j}b_i=&
\sum_{i=1}^{m+z-r-x}  v_i-\sum_{i=1}^{m-r}  u_i+\sum_{i=1}^{  r+x}f_i\\&-\sum_{i=1}^{  r+x}\max\{e_{i-x-j}, f_i\}=
\sum_{i=1}^{j}\hat b_i.
\end{array}
$$
Therefore $\bb=\hat \bb$ and (\ref{eqrmimajpolhatsA}) is equivalent to (\ref{eqrmimajpol}).
\end{itemize}
\hfill $\Box$

\medskip

In the following example we show that, in general,  conditions (\ref{eqetapol}) and (\ref{eqinteriedpol})-(\ref{eqrmimajpolhatsA})
  are not sufficient if the field is not algebraically closed and $A>0$.

\begin{example}
  \label{exalgclosed}
  Let $P(s)=
  \begin{bmatrix}s&1\\-1&s
  \end{bmatrix}
  \in \FF[s]^{2\times2}
  $. Notice that
  $$
  d=1, \quad r=2, \quad \bc=\emptyset , \quad \bu=\emptyset, \quad \phi_1(s,t)=1, \quad \phi_2(s,t)=s^2+t^2,
  $$$$e_1=e_2=0, \quad \alpha_1(s)=1, \quad \alpha_2(s)=s^2+1.
  $$ 
  Let $x=0$, $z=1$. We prescribe $$f_1=f_2=0, \quad \bd=\emptyset, \quad \bv=(1).$$
  Then, $A=1$, $\hat b_1=1$ and
   conditions (\ref{eqetapol}) and (\ref{eqinteriedpol})-(\ref{eqrmimajpolhatsA}) hold.
  Hence,  there exists a pencil $W(s)\in \efe[s]^{1\times 2}$ such that   $\begin{bmatrix}P(s)\\W(s)\end{bmatrix}$ has $v_1=1$ as row minimal index, it does not have column minimal indices, and $f_1=f_2=0$ are its partial multiplicities of $\infty$. As
$\rank \left(\begin{bmatrix}P(s)\\W(s)\end{bmatrix}\right) =2$, if
  $\beta_1(s) \mid \beta_2(s)$ are its invariant factors then, by Lemma \ref{theoDeDoMa},
  $\deg(\beta_1)+ \deg(\beta_2)=1$, and by Theorem \ref{theoprescr4}, $\beta_1(s) \mid \alpha_1(s)=1$ and
  $\beta_2(s) \mid \alpha_2(s)=s^2+1$, i.e., there exists a
   polynomial $\beta_2(s)$ such that $\beta_2(s) \mid
  s^2+1$ and $\deg(\beta_2)=1$.

  If $\efe=\RR$ there is no such polynomial. Therefore, conditions (\ref{eqetapol}) and (\ref{eqinteriedpol})-(\ref{eqrmimajpolhatsA})
  are not sufficient if the field is not algebraically closed and $A>0$.
  
  If $\efe =\CC$, there exist  only two possible chains of invariant factors, $\beta_1(s)=1$, 
  $\beta_2(s)=s+i$, or $\beta_1(s)=1$, $\beta_2(s)=s-i$, that admit  completion  having the prescribed invariants.  For example,
  $$
  \begin{bmatrix}P(s)\\W(s)\end{bmatrix}=
  \begin{bmatrix}s&1\\-1&s\\0&s\pm  i
  \end{bmatrix}
  \in \CC[s]^{3\times2},
  $$ has
  $v_1=1$ as row minimal index, it does not have column minimal indices, and $f_1=f_2=0$ are its partial multiplicities of $\infty$.
  \end{example}
 
In the following  two corollaries, we solve Problem \ref{problem} when the prescribed invariants are  the infinite elementary divisors and the column minimal indices, and then  the infinite elementary divisors and the row minimal indices.
We remark that Corollary \ref{corprescrdeicmi} is valid on arbitrary fields, while the sufficiency part of
Corollary \ref{corprescrdeirmi}
requires  the field to  be algebraically closed.

\begin{corollary} {\rm (Prescription of infinite structure and column minimal indices)}
  \label{corprescrdeicmi}
  Let $P(s)\in\efe[s]^{m\times n}$, $\deg(P(s))=d$, $\rank (P(s))=r$.
  Let  $e_1\leq\dots\leq e_r$ be the partial multiplicities of $\infty$ and
  $\bc=(c_1,  \dots,  c_{n-r})$  the  column minimal indices of $P(s)$.

Let $z$ and $x$ be integers such that $0\leq x\leq \min\{z, n-r\}$ and 
let  $f_1\leq\dots \leq f_{r+x} $ and
 $d_1 \geq  \dots  \geq  d_{n-r-x}\geq 0$ 
  be non negative integers.
  There exists a polynomial matrix  $W(s)\in \efe[s]^{z\times n}$ such that
$\deg(W(s))\leq d$, $\rank\left(\begin{bmatrix}P(s)\\W(s)\end{bmatrix}\right)=r+x$ 
and $\begin{bmatrix}P(s)\\W(s)\end{bmatrix}$ has
 $f_1\leq\dots \leq f_{r+x} $ as partial multiplicities of $\infty$ and 
 $\bd=(d_1,\dots, d_{n-r-x})$ as column minimal indices
if and only if
(\ref{eqinteriedpol}),
\begin{equation}\label{eqcde}
 \sum_{i=1}^{n-r}c_i-\sum_{i=1}^{n-r-x}d_i\geq \sum_{i=1}^{r+x}\max\{e_{i-x}, f_i\}-\sum_{i=1}^{r}e_i-xd,
\end{equation}
\begin{equation}\label{eqcmimajtildea}  \bc \prec'  (\bd , \tilde \ba),\end{equation}
where $\tilde \ba = (\tilde a_1, \dots,  \tilde a_x )$ is
\begin{equation}\label{eqaiedcmi} 
\begin{array}{rl}\tilde a_1=&
\sum_{i=1}^{n-r}c_i-\sum_{i=1}^{n-r-x}d_i\\&+\sum_{i=1}^{r}e_i-\sum_{i=1}^{r+x-1}\max\{e_{i-x+1}, f_i\}+(x-1)d,\\ \tilde a_j=&
\sum_{i=1}^{r+x-j+1}\max\{e_{i-x+j-1}, f_i\}-\sum_{i=1}^{r+x-j}\max\{e_{i-x+j}, f_i\}
  -d,\\& \hfill  2\leq j \leq x.
\end{array}
\end{equation}
\end{corollary}
{\bf Proof.}
First notice that by (\ref{eqcde}) and Remark \ref{remposttheo}.\ref{remitdecreasing}, $\tilde \ba$ is 
well defined.

Let $\bu=(u_1, \dots, u_{m-r})$ with
$u_1 \geq  \dots  \geq  u_\eta>0=u_{\eta+1}=\dots=u_{m-r}$ be the row minimal indices of $P(s)$.

Assume that  there exists $W(s)\in \efe[s]^{z\times n}$ 
such that
$\deg(W(s))\leq d$, $\rank\left(\begin{bmatrix}P(s)\\W(s)\end{bmatrix}\right)=r+x$ 
and $\begin{bmatrix}P(s)\\W(s)\end{bmatrix}$ has
 $f_1\leq\dots \leq f_{r+x} $ as partial multiplicities of $\infty$ and 
 $\bd=(d_1,\dots, d_{n-r-x})$ as column minimal indices. Let
$\bv=(v_1,\dots, v_{m+z-r-x})$  be its row minimal indices. By Theorem \ref{theoprescrdeirmicmi} and Remark \ref{remacdecr}.\ref{remitac}, 
(\ref{eqinteriedpol}), (\ref{eqvusA}) and (\ref{eqcmimajpolhatsA}) hold, where
$A$ and  $\hat \ba$ are defined in (\ref{eqbigAdeicmirmi}) and  (\ref{eqadeicmirmi}), respectively. Observe that
$\hat a_j=\tilde a_j$ for $2\leq j \leq x$.
 As 
$$
\sum_{i=1}^{m+z-r-x}v_i-\sum_{i=1}^{m-x}u_i+ \sum_{i=1}^{r+x}f_i=A+\sum_{i=1}^{n-r}c_i-\sum_{i=1}^{n-r-x}d_i+\sum_{i=1}^{r}e_i+xd,
$$
we have  $\hat \ba=\tilde \ba$, and from  (\ref{eqvusA}) and  (\ref{eqcmimajpolhatsA}) we obtain (\ref{eqcde}) and  (\ref{eqcmimajtildea}), respectively.

Conversely, assume that (\ref{eqinteriedpol}),
(\ref{eqcde}) and 
(\ref{eqcmimajtildea}) hold.
Define
$\tilde \bb = (\tilde b_1, \dots, \tilde b_{z-x} )$ as
$$\begin{array}{rl}\tilde b_1=&
{ \sum_{i=1}^{n-r}c_i-\sum_{i=1}^{n-r-x}d_i+\sum_{i=1}^{r}e_i+xd
-\sum_{i=1}^{  r+x}\max\{e_{i-x-1}, f_i\},}
  \\ \tilde b_j=&{ 
  \sum_{i=1}^{  r+x}\max\{e_{i-x-j+1}, f_i\}-\sum_{i=1}^{  r+x}\max\{e_{i-x-j}, f_i\}
  ,\quad 2\leq j \leq z-x.}
\end{array}
$$
By (\ref{eqcde}) and Remark \ref{remposttheo}.\ref{remitdecreasing}, $\tilde \bb$ is 
well defined.
Take
$\bv=\bu \cup \tilde \bb$ and $\bar \eta= \#\{i\;  : \; v_i>0\}$. Then (\ref{eqetapol}) holds.
Moreover,
$$
\begin{array}{rrl}
\sum_{i=1}^{m+z-r-x}v_i-\sum_{i=1}^{m-r}u_i=&\sum_{i=1}^{z-x}\tilde b_i=&\sum_{i=1}^{n-r}c_i-\sum_{i=1}^{n-r-x}d_i+\sum_{i=1}^{r}e_i
\\&&-\sum_{i=1}^{  r+x}
\max\{e_{i-z}, f_i\}+xd.
\end{array}
$$
From (\ref{eqinteriedpol}) we obtain
\begin{equation}\label{sumvmu}
\sum_{i=1}^{m+z-r-x}v_i-\sum_{i=1}^{m-r}u_i=\sum_{i=1}^{n-r}c_i-\sum_{i=1}^{n-r-x}d_i+\sum_{i=1}^{r}e_i-\sum_{i=1}^{  r+x}f_i+xd.
\end{equation}
Let $\alpha_1(s)\mid\cdots\mid\alpha_r(s)$ be the
  invariant factors of $P(s)$ and 
let $A$, $\hat \ba$ and  $\hat \bb$ be as in (\ref{eqbigAdeicmirmi}), (\ref{eqadeicmirmi}) and (\ref{eqbdeicmirmi}), respectively. Then $\hat a_j=\tilde a_j$ for $2\leq j \leq x$.
From (\ref{sumvmu}) we obtain $A=0$; hence (\ref{eqAleq}) holds, and from (\ref{eqcde}) we derive (\ref{eqvusA}). 
From (\ref{sumvmu}) we get
$$
  \hat a_1=\sum_{i=1}^{n-r}c_i-\sum_{i=1}^{n-r-x}d_i+\sum_{i=1}^{r}e_i-\sum_{i=1}^{r+x-1}\max\{e_{i-x+1}, f_i\}+(x-1)d
  =\tilde a_1;
$$
hence (\ref{eqcmimajtildea}) and  (\ref{eqcmimajpolhatsA}) are equivalent.
We have 
$$
\begin{array}{rl}
  \sum_{i=1}^j \hat b_i=&\sum_{i=1}^{m+z-r-x}v_i-\sum_{i=1}^{m-r}u_i+\sum_{i=1}^{  r+x}f_i-\sum_{i=1}^{  r+x}\max\{e_{i-x-j}, f_i\}\\
  =&\sum_{i=1}^{n-r}c_i-\sum_{i=1}^{n-r-x}d_i+\sum_{i=1}^{r}e_i+xd\\&-\sum_{i=1}^{  r+x}\max\{e_{i-x-j}, f_i\}
  =\sum_{i=1}^j \tilde b_i, \quad 1\leq j \leq z-x.
\end{array}
$$
Therefore $\hat \bb=\tilde \bb$; hence  $\bv=\bu \cup \hat \bb$. By Lemma \ref{lemmacup}, (\ref{eqrmimajpolhatsA}) holds.
By Theorem \ref{theoprescrdeirmicmi} and Remark \ref{remacdecr}.\ref{remitac}, the sufficiency of the conditions (\ref{eqinteriedpol}),
(\ref{eqcde}) and 
(\ref{eqcmimajtildea}) follows.
\hfill $\Box$

\begin{example}
  \label{exiedcmi}
 Let $P(s)=
  \begin{bmatrix}s&1\\-1&s
  \end{bmatrix}
  \in \FF[s]^{2\times2}
  $ as in Example \ref{exalgclosed}, with 
  $
  d=1$, $r=2$, $e_1=e_2=0$, $\bc=\emptyset$.
  Let $x=0$, $z=1$. We prescribe $f_1=f_2=0$ and $\bd=\emptyset$.
  Conditions (\ref{eqinteriedpol}),
(\ref{eqcde}) and 
 (\ref{eqcmimajtildea}) hold.
  Hence,  there exists a pencil $W(s)\in \efe[s]^{1\times 2}$ such that   $\begin{bmatrix}P(s)\\W(s)\end{bmatrix}$ has  $f_1=f_2=0$ as partial multiplicities of $\infty$ and it does not have column minimal indices. 
  
  As
$\rank \left(\begin{bmatrix}P(s)\\W(s)\end{bmatrix}\right) =2$, if
  $\beta_1(s) \mid \beta_2(s)$ are its invariant factors and $v_1$ its row minimal index, then, by Lemma \ref{theoDeDoMa},
  $\deg(\beta_1)+ \deg(\beta_2)+v_1=2$, and by Theorem \ref{theoprescr4}, $\beta_1(s) \mid \alpha_1(s)=1$ and
$\beta_2(s) \mid \alpha_2(s)=s^2+1$.
If $W(s)=0\in  \efe[s]^{1\times 2}$ we obtain $v_1=0$, $\beta_2(s)=s^2+1$. If $W(s)=
\begin{bmatrix}1&0\end{bmatrix}
\in  \efe[s]^{1\times 2}$ we obtain $v_1=2$, $\beta_2(s)=1$. In both cases the completion has the prescribed invariants.
As we have seen in  Example \ref{exalgclosed}, if $\efe$ is an algebraically closed field, we can also obtain  $v_1=1$ and $\deg(\beta_2)=1$. 
  \end{example}

\begin{corollary}{\rm (Prescription of infinite structure and row minimal indices)}
  \label{corprescrdeirmi}
  Let $\efe$ be an algebraically closed field.
Let $P(s)\in\efe[s]^{m\times n}$,  $\deg(P(s))=d$, $\rank (P(s))=r$.
  Let $\alpha_1(s)\mid\cdots\mid\alpha_r(s)$ be the invariant factors,
$e_1\leq\cdots\leq e_r$
the partial multiplicities of $\infty$, 
$\bc=(c_1,  \dots,  c_{n-r})$  the
column minimal indices, and $ \bu=(u_1, \dots,  u_{m-r})$  the row minimal
indices of $P(s)$, where 
$u_1 \geq  \dots  \geq  u_{\eta}  > u_{\eta +1} = \dots =u_{m-r} = 0$.

Let $z$ and $x$ be integers such that $0\leq x\leq \min\{z, n-r\}$ and 
let  $f_1\leq\dots \leq f_{r+x} $ and
$v_1 \geq  \dots  \geq  v_{\bar \eta} > v_{\bar \eta+1}=\dots= v_{m+z-r-x}=0$ 
  be non negative integers. Let
\begin{equation}\label{eqbigtildeAiedrmi}
  \tilde A=\sum_{i=1}^{r+x}f_i-\sum_{i=1}^{r}e_i-\sum_{i=1}^{x}c_{i}+\sum_{i=1}^{m+z-r-x}v_i-\sum_{i=1}^{m-r}u_i-xd.
  \end{equation}
There exists $W(s)\in \efe[s]^{z\times n}$ such that
$\deg(W(s))\leq d$, $\rank \left(\begin{bmatrix}P(s)\\W(s)\end{bmatrix} \right)=r+x$
and $\begin{bmatrix}P(s)\\W(s)\end{bmatrix}$  has 
 $f_1\leq\dots \leq f_{r+x} $ as  partial multiplicities of $\infty$ 
 and
$\bv=(v_1,\dots, v_{m+z-r-x})$ as row minimal indices
if and only if 
 (\ref{eqetapol}), (\ref{eqinteriedpol}),
\begin{equation}\label{eqtildeAleq}
\tilde A\leq \sum_{i=r+x-z+1}^{r}\deg(\alpha_i),
\end{equation}
\begin{equation}\label{eqvustildeA}
 \sum_{i=1}^{m+z-r-x}v_i-\sum_{i=1}^{m-r} u_i\geq \max\{0, \tilde A\}+\sum_{i=1}^{r+x}\max\{e_{i-x}, f_i\}-\sum_{i=1}^{r+x}f_i,
\end{equation}
\begin{equation}\label{eqcmimajpoltildeA}  (c_1, \dots, c_x) \prec   \tilde  \ba,\end{equation}
\begin{equation}\label{eqrmimajpoltildeA} \bv \prec'  (\bu , \tilde \bb),\end{equation}
where $\tilde  \ba = (\tilde a_1, \dots, \tilde a_x )$ and $\tilde \bb = (\tilde b_1, \dots,  \tilde b_{z-x})$ 
are 
\begin{equation}\label{eqtildeaiedrmi}
\begin{array}{rl}\tilde  a_1=&\sum_{i=1}^{m+z-r-x}v_i-\sum_{i=1}^{m-r}u_i\\&+\sum_{i=1}^{r+x}f_i-\sum_{i=1}^{r+x-1}\max\{e_{i-x+1}, f_i\}-\tilde A-d,\\ \tilde  a_j=&
\sum_{i=1}^{r+x-j+1}\max\{e_{i-x+j-1}, f_i\}\\&-\sum_{i=1}^{r+x-j}\max\{e_{i-x+j}, f_i\}
  -d,\quad 2\leq j \leq x,
\end{array}
\end{equation}
\begin{equation}\label{eqtildebiedrmi}
\begin{array}{rl}\tilde  b_1=&   {\sum_{i=1}^{m+z-r-x}v_i-\sum_{i=1}^{m-r}u_i+\min\{0,\deg(\alpha_{r})-\tilde A \}}\\& {+\sum_{i=1}^{  r+x}f_i-\sum_{i=1}^{  r+x}\max\{e_{i-x-1}, f_i\},}
  \\ \tilde b_j=&  {\min\{0,\sum_{i=r-j+1}^{r}\deg(\alpha_{i})-\tilde A\}-\min\{0,\sum_{i=r-j+2}^{r}\deg(\alpha_{i})-\tilde A \}}\\&
  { +\sum_{i=1}^{  r+x}\max\{e_{i-x-j+1}, f_i\}-\sum_{i=1}^{  r+x}\max\{e_{i-x-j}, f_i\}
  ,\quad 2\leq j \leq z-x.}
\end{array}
\end{equation}
\end{corollary}
\begin{remark}\label{remac2}
As in Theorem \ref{theoprescrdeirmicmi}, the necessity of the conditions, and the sufficiency when $\tilde A\leq0$ hold for arbitrary fields.
  \end{remark}   
{\bf Proof of Corollary \ref{corprescrdeirmi}.} 
Assume that there exists a polynomial matrix  $W(s)\in \efe[s]^{z\times n}$  such that
$\deg(W(s))\leq d$, $\rank \left(\begin{bmatrix}P(s)\\W(s)\end{bmatrix}\right) =r+x$
and $\begin{bmatrix}P(s)\\W(s)\end{bmatrix}$  has 
 $f_1\leq\dots \leq f_{r+x} $ as  partial multiplicities of $\infty$, 
 and
$\bv=(v_1,\dots, v_{m+z-r-x})$ as row minimal indices. Let
$\bd=(d_1,\dots, d_{n-r-x})$ be its column minimal indices. By Theorem \ref{theoprescrdeirmicmi}, 
(\ref{eqetapol}) and (\ref{eqinteriedpol})-(\ref{eqrmimajpolhatsA}) hold,
where
$A$, $\hat \ba$ and $\hat \bb$ are defined in (\ref{eqbigAdeicmirmi}), (\ref{eqadeicmirmi}) and (\ref{eqbdeicmirmi}), respectively.
From (\ref{eqcmimajpolhatsA}) we have $d_i\geq c_{i+x}$, $1\leq i \leq n-r-x$; hence $\tilde A=A+\sum_{i=1}^{n-r-x}c_{i+x}-\sum_{i=1}^{n-r-x}d_i\leq A$.
From (\ref{eqAleq}) and (\ref{eqvusA})  we obtain (\ref{eqtildeAleq}) and (\ref{eqvustildeA}), respectively.
Moreover,
$$\begin{array}{rl}
  \sum_{i=1}^j\hat b_i=&\sum_{i=1}^{m+z-r-x}v_i-\sum_{i=1}^{m-r}u_i+\sum_{i=1}^{  r+x}f_i-\sum_{i=1}^{  r+x}\max\{e_{i-x-j}, f_i\}\\&+\min\{0,\sum_{i=r-j+1}^{r}\deg(\alpha_{i})-A\}
  \\
  \leq &\sum_{i=1}^{m+z-r-x}v_i-\sum_{i=1}^{m-r}u_i+\sum_{i=1}^{  r+x}f_i-\sum_{i=1}^{  r+x}\max\{e_{i-x-j}, f_i\}\\&+\min\{0,\sum_{i=r-j+1}^{r}\deg(\alpha_{i})-\tilde A\}=\sum_{i=1}^j\tilde b_i,\quad 1\leq j \leq z-x,\\
  \sum_{i=1}^{z-x}\hat b_i=&\sum_{i=1}^{m+z-r-x}v_i-\sum_{i=1}^{m-r}u_i=\sum_{i=1}^{z-x}\tilde b_i,
\end{array}
$$
i.e., $\hat \bb\prec \tilde \bb$. From (\ref{eqrmimajpolhatsA}) 
we obtain (\ref{eqrmimajpoltildeA}) (see  Remark \ref{gmajcp}).

We have
$$
\sum_{i=1}^x \tilde a_i =\sum_{i=1}^{m+z-r-x}v_i-\sum_{i=1}^{m-x}u_i+\sum_{i=1}^{r+x}f_i-\sum_{i=1}^{r}e_i-\tilde A-xd=\sum_{i=1}^{x}c_i.
$$
Let $j\in \{1, \dots, x\}$.
Thus,
$$
\begin{array}{rl}
  \sum_{i=1}^j \hat a_i=&\sum_{i=1}^{m+z-r-x}v_i-\sum_{i=1}^{m-r}u_i+\sum_{i=1}^{  r+x}f_i-A-jd\\&-\sum_{i=1}^{r+x-j}\max\{e_{i-x+j}f_i\}\\=&
  \sum_{i=1}^j \tilde a_i+\tilde A-A=
\sum_{i=1}^j \tilde a_i+\sum_{i=1}^{n-r-x}c_{i+x}-\sum_{i=1}^{n-r-x}d_i.
\end{array}
$$
Let 
$\hat h_j=\min\{i\; : \; d_{i-j+1}<c_i\}$. Then   $j \leq \hat h_j \leq n-r-x+j$. From (\ref{eqcmimajpolhatsA}),
$$
\begin{array}{rl}
  \sum_{i=1}^{\hat h_j}c_i\leq &\sum_{i=1}^{\hat h_j-j}d_i+ \sum_{i=1}^j \hat a_i=\sum_{i=1}^j \tilde a_i+\sum_{i=1}^{n-r-x}c_{i+x}-\sum_{i=\hat h_j-j+1}^{n-r-x}d_i\\
  = & \sum_{i=1}^j \tilde a_i+\sum_{i=1}^{\hat h_j-j}c_{i+x}+\sum_{i=\hat h_j-j+1}^{n-r-x}(c_{i+x}-d_i)\\\leq& \sum_{i=1}^j \tilde a_i+\sum_{i=1}^{\hat h_j-j}c_{i+x};
\end{array}
$$
hence
$$\begin{array}{rl}
  \sum_{i=1}^{j}c_i\leq &\sum_{i=1}^j \tilde a_i+\sum_{i=1}^{\hat h_j-j}c_{i+x}-\sum_{i=j+1}^{\hat h_j}c_{i}\\=& \sum_{i=1}^j \tilde a_i+\sum_{i=1}^{\hat h_j-j}(c_{i+x}-c_{i+j})\leq \sum_{i=1}^j \tilde a_i.
\end{array}
$$
Therefore, (\ref{eqcmimajpoltildeA}) holds.

Conversely, let us assume that
(\ref{eqetapol}),
(\ref{eqinteriedpol}), and 
(\ref{eqtildeAleq})-(\ref{eqrmimajpoltildeA}) hold, where $\tilde A$, $\tilde \ba$ and   $\tilde \bb$ are defined in (\ref{eqbigtildeAiedrmi}), (\ref{eqtildeaiedrmi}) and
 (\ref{eqtildebiedrmi}), respectively. 
Define $d_i=c_{i+x}$, $1\leq i\leq n-r-x$, and $\bd=(d_1, \dots, d_{n-r-x})$.
Let $A$, $\hat \ba$ and $\hat \bb$ be as in 
(\ref{eqbigAdeicmirmi}), (\ref{eqadeicmirmi}) and (\ref{eqbdeicmirmi}), respectively.
Then $A=\tilde A$, $\hat \ba=\tilde \ba$  and $\hat \bb=\tilde \bb$; hence  (\ref{eqtildeAleq}),
(\ref{eqvustildeA}) and (\ref{eqrmimajpoltildeA}) are equivalent to  (\ref{eqAleq}),
(\ref{eqvusA}) and (\ref{eqrmimajpolhatsA}), respectively.
We have $d_i\geq c_{i+x}$ for $1\leq i \leq n-r-x$ and $\sum_{i=1}^{x}\hat a_i=
\sum_{i=1}^{n-r}c_i - \sum_{i=1}^{n-r-x}d_i$.

Let $j\in \{1, \dots, x\}$ and 
$\hat h_j=\min\{i\; : \; d_{i-j+1}<c_i\}$. 
Then $j \leq \hat h_j \leq n-r-x+j$. 
If
 $j\leq i \leq \hat h_j-1$ then 
$c_{i+x-j+1}=d_{i-j+1}\geq c_i$; hence
$
c_{i+x-j+1}=c_i$, $j\leq i \leq \hat h_j-1.
$
Thus, 
$$
\begin{array}{rl}
\sum_{i=1}^{\hat h_j}c_i=&\sum_{i=1}^{j-1}c_i+\sum_{i=j}^{\hat h_j-1}c_{i+x-j+1}+c_{\hat h_j}=
\sum_{i=1}^{j-1}c_i+c_{\hat h_j}+\sum_{i=1}^{\hat h_j-j}c_{i+x}\\=&
\sum_{i=1}^{j-1}c_i+c_{\hat h_j}+\sum_{i=1}^{\hat h_j-j}d_{i}\leq \sum_{i=1}^{j}c_i+\sum_{i=1}^{\hat h_j-j}d_{i}.
\end{array}
$$
From (\ref{eqcmimajpoltildeA}), $\sum_{i=1}^{j}c_i\leq \sum_{i=1}^{j}\tilde a_i$; hence 
$
\sum_{i=1}^{\hat h_j}c_i
\leq \sum_{i=1}^{\hat h_j-j}d_{i}+\sum_{i=1}^{j}\tilde a_i.
$
Therefore, (\ref{eqcmimajpolhatsA}) holds.
By Theorem \ref{theoprescrdeirmicmi}, the sufficiency of conditions (\ref{eqetapol}),
(\ref{eqinteriedpol}), and 
(\ref{eqtildeAleq})-(\ref{eqrmimajpoltildeA})  follows.
\hfill $\Box$

\medskip

Example  \ref{exalgclosed}  also shows that, in general,  conditions (\ref{eqetapol}),
(\ref{eqinteriedpol}), and 
(\ref{eqtildeAleq})-(\ref{eqrmimajpoltildeA}) 
  are not sufficient if the field is not algebraically closed and $\tilde A>0$.
  
  \subsection{Prescription of finite and singular structures}
  
In this subsection we prescribe first the invariant factors and the row and column minimal indices, then  the invariant factors and the column minimal indices, and afterwards  the invariant factors and the row  minimal indices.
  The results are valid over arbitrary fields.
  
\begin{theorem} {\rm (Prescription of finite and singular structures)}
\label{theoprescrfifrmicmi}
\hspace{-0.07cm}Let $P(s)\hspace{-0.06cm}\in\hspace{-0.06cm}\efe[s]^{m\times n}$,  $\deg(P(s))=d$, $\rank (P(s))=r$.
  Let $\alpha_1(s)\mid\cdots\mid\alpha_r(s)$ be the invariant factors,
$e_1\leq\cdots\leq e_r$
the partial multiplicities of $\infty$, 
$\bc=(c_1,  \dots,  c_{n-r})$  the
column minimal indices, and $ \bu=(u_1, \dots,  u_{m-r})$  the row minimal
indices of $P(s)$, where 
$u_1 \geq  \dots  \geq  u_{\eta}  > u_{\eta +1} = \dots =u_{m-r} = 0$.

Let $z$ and $x$ be integers such that $0\leq x\leq \min\{z, n-r\}$ and 
let  $\beta_1(s)\mid\dots \mid \beta_{r+x}(s) $ be monic polynomials and
 $d_1 \geq  \dots  \geq  d_{n-r-x}\geq 0$ and
 $v_1 \geq  \dots  \geq  v_{\bar \eta} > v_{\bar \eta+1}=\dots= v_{m+z-r-x}=0$  non negative integers.
Let
\begin{equation}\label{eqbigBfifcmirmi}
B=\sum_{i=1}^{r+x}\deg(\beta_i)-\sum_{i=1}^{r}\deg(\alpha_i)+\sum_{i=1}^{n-r-x}d_i-\sum_{i=1}^{n-r}c_i+\sum_{i=1}^{m+z-r-x}v_i-\sum_{i=1}^{m-r}u_i-xd
.\end{equation}
There exists $W(s)\in \efe[s]^{z\times n}$ such that $\deg(W(s))\leq d$, $\rank \left(\begin{bmatrix}P(s)\\W(s)\end{bmatrix}\right) =r+x$
and $\begin{bmatrix}P(s)\\W(s)\end{bmatrix}$  has   
 $\beta_1(s)\mid\dots \mid \beta_{r+x}(s) $  as  invariant factors,
 $\bd=(d_1,\dots, d_{n-r-x})$ as column minimal indices and 
$\bv=(v_1,\dots, v_{m+z-r-x})$ as row minimal indices
if and only if 
(\ref{eqetapol}), (\ref{eqinterfifpol}),
\begin{equation}\label{eqBleq}
B\leq \sum_{i=r+x-z+1}^{r}e_i,
\end{equation}
\begin{equation}\label{eqvusB}
 \sum_{i=1}^{m+z-r-x}v_i-\sum_{i=1}^{m-r} u_i\geq \max\{0, B\}+\sum_{i=1}^{r+x}\deg(\lcm(\alpha_{i-x}, \beta_i))-\sum_{i=1}^{r+x}\deg(\beta_i),
\end{equation}
\begin{equation}\label{eqcmimajpolhatsB}  \bc \prec'  (\bd , \hat \ba),\end{equation}
\begin{equation}\label{eqrmimajpolhatsB} \bv \prec'  (\bu , \hat \bb),\end{equation}
where $\hat \ba = (\hat a_1, \dots, \hat a_x )$ and $\hat \bb = (\hat b_1, \dots,  \hat b_{z-x})$ 
are 
\begin{equation}\label{eqafifcmirmi}
\begin{array}{rl}\hat a_1=&  \sum_{i=1}^{m+z-r-x}v_i-\sum_{i=1}^{m-r}u_i+\sum_{i=1}^{r+x}\deg(\beta_i)\\&-\sum_{i=1}^{r+x-1}\deg(\lcm(\alpha_{i-x+1}, \beta_i))-B-d
  ,\\ \hat a_j=& 
\sum_{i=1}^{r+x-j+1}\deg(\lcm(\alpha_{i-x+j-1}, \beta_i))\\&-\sum_{i=1}^{r+x-j}\deg(\lcm(\alpha_{i-x+j}, \beta_i))
  -d,\quad 2\leq j \leq x,
\end{array}
\end{equation}
\begin{equation}\label{eqbfifcmirmi}
\begin{array}{rl}\hat b_1=&  {\sum_{i=1}^{m+z-r-x}v_i-\sum_{i=1}^{m-r}u_i+\min\{0,e_r-B \}}\\& {+\sum_{i=1}^{  r+x}\deg(\beta_i)-\sum_{i=1}^{  r+x}\deg(\lcm(\alpha_{i-x-1}, \beta_i)),}
  \\ \hat b_j=& {\min\{0,\sum_{i=r-j+1}^{r}e_i-B\}-\min\{0,\sum_{i=r-j+2}^{r}e_i-B \}}\\& 
  +\sum_{i=1}^{  r+x}\deg(\lcm(\alpha_{i-x-j+1}, \beta_i))-\sum_{i=1}^{  r+x}\deg(\lcm(\alpha_{i-x-j}, \beta_i)),\\& \hfill
   2\leq j \leq z-x.
\end{array}
\end{equation}
\end{theorem}
{\bf Proof.}
The proof is analogous to that of Theorem  \ref{theoprescrdeirmicmi} by exchanging the roles of the invariant factors and the infinite elementary divisors. 
In the sufficiency part,
we take $g=\min\{k\geq 0\; : \;B\leq \sum_{i=r-k+1}^{r} e_i \}$.
\begin{itemize}
  \item  
    If  $B>0$, then
$0<B-\sum_{i=r-g+2}^{r}e_i\leq e_{r-g+1}
    .$
    
Let $h=\min\{k\; : \;B-\sum_{i=r-g+2}^{r}e_i\leq e_{k-g+1}\}$ and
$w=e_{h-g}+e_{h-g+1}+\sum_{i=r-g+2}^{r}e_i-B.$ Then $e_{h-g}\leq w<e_{h-g+1}$.
Define 
$$
\begin{array}{rll}
  f_i=&e_{i-x-g},&1\leq i\leq h+x-1,
  \\
  f_{h+x}=&w,\\
  f_i=&e_{i-x-g+1},&h+x+1\leq i\leq r+x.
\end{array}
$$
Observe that $f_1, \dots, f_{r+x}$  can be defined in an arbitrary field.
\item
 If  $B\leq 0$, 
  define
$$
\begin{array}{rll}
  f_i=&e_{i-x},&1\leq i\leq r+x-1,
  \\
  f_{r+x}=&e_r-B.
\end{array}
$$
\end{itemize}

In both cases, if $
\gamma_i(s,t)=t^{f_i}t^{\deg(\beta_i)}\beta_i(\frac{s}{t}) 
$, $1\leq i \leq r+x$, then
$\gamma_1(s,t)\mid \dots \mid \gamma_{r+x}(s,t)$
and, arguing as in the proof of Theorem  \ref{theoprescrdeirmicmi}, 
(\ref{eqinterfipolhom})-(\ref{eqdegsumpol}) are satisfied.
By Theorem \ref{theoprescr4}, the sufficiency of conditions (\ref{eqetapol}), (\ref{eqinterfifpol}),
(\ref{eqBleq})-(\ref{eqrmimajpolhatsB}) follows.
\hfill $\Box$

\begin{example}
  \label{exprifcmirmi}
  Let $P(s)=
  \begin{bmatrix}s^2&-1&0\\0&0&0
  \end{bmatrix}
  \in \FF[s]^{2\times3}
  $. Then
  $$
  d=2, \quad r=1, \quad \alpha_1(s)=\alpha_2(s)=1, \quad e_1=e_2=0, \quad \bc=(2,0) , \quad \bu=(0).$$
  Let $x=z=1$. 
  In this case, conditions (\ref{eqetapol}), (\ref{eqinterfifpol}),
(\ref{eqBleq}), (\ref{eqvusB}) and (\ref{eqrmimajpolhatsB}) reduce to $\beta_1(s)=1$, $v_1=u_1=0$, $\deg(\beta_1)+d_1\leq 4$, and condition (\ref{eqcmimajpolhatsB}) is equivalent to $d_1\geq 2$ or $d_1=0$.

For example, if we prescribe $\beta_1(s)=\beta_2(s)=1$, $v_1=0$ and $d_1=3$, a possible completion is
$Q(s)=\begin{bmatrix}s^2&-1&0\\0&0&0\\0&s&-1
  \end{bmatrix}$.
  \end{example}

\begin{corollary} {\rm (Prescription of finite structure and column minimal indices)}
  \label{corprescrfifcmi}
Let $P(s)\in\efe[s]^{m\times n}$,  $\deg(P(s))=d$, $\rank (P(s))=r$.
  Let $\alpha_1(s)\mid\cdots\mid\alpha_r(s)$ be the invariant factors,
and
$\bc=(c_1,  \dots,  c_{n-r})$  the
column minimal indices of $P(s)$.

Let $z$ and $x$ be integers such that $0\leq x\leq \min\{z, n-r\}$ and 
let  $\beta_1(s)\mid \dots \mid \beta_{r+x}(s) $ be monic polynomials  and
 $d_1 \geq  \dots  \geq  d_{n-r-x}\geq 0$ 
   non negative integers.
There exists a polynomial matrix  $W(s)\in \efe[s]^{z\times n}$ such that $\deg(W(s))\leq d$, $\rank\left(\begin{bmatrix}P(s)\\W(s)\end{bmatrix}\right) =r+x$
and $\begin{bmatrix}P(s)\\W(s)\end{bmatrix}$  has
 $\beta_1(s)\mid \dots \mid \beta_{r+x}(s) $ as invariant factors  and 
 $\bd=(d_1,\dots, d_{n-r-x})$ as column minimal indices
if and only if (\ref{eqinterfifpol}),
\begin{equation}\label{eqcdefif}
 \sum_{i=1}^{n-r}c_i-\sum_{i=1}^{n-r-x}d_i\geq \sum_{i=1}^{r+x}\deg(\lcm(\alpha_{i-x}, \beta_i))-\sum_{i=1}^{r}\deg(\alpha_i)-xd,
\end{equation}
\begin{equation}\label{eqcmimajtildeafif}  \bc \prec'  (\bd , \tilde \ba),\end{equation}
where $\tilde \ba = (\tilde a_1,  \dots,  \tilde a_x )$ is
\begin{equation}\label{eqafifcmi}
\begin{array}{rl}\tilde a_1=& 
\sum_{i=1}^{n-r}c_i-\sum_{i=1}^{n-r-x}d_i\\&+\sum_{i=1}^{r}\deg(\alpha_i)-\sum_{i=1}^{r+x-1}\deg(\lcm(\alpha_{i-x+1}, \beta_i))+(x-1)d
  ,\\ \tilde a_j=& 
\sum_{i=1}^{r+x-j+1}\deg(\lcm(\alpha_{i-x+j-1}, \beta_i))\\&-\sum_{i=1}^{r+x-j}\deg(\lcm(\alpha_{i-x+j}, \beta_i))
  -d,\quad 2\leq j \leq x.
\end{array}
\end{equation}
\end{corollary}
{\bf Proof.} The result is a consequence of Theorem \ref{theoprescrfifrmicmi}. The proof follows the scheme of Corollary \ref{corprescrdeicmi} by exchanging the roles of the invariant factors and the infinite elementary divisors. 
\hfill $\Box$

 \begin{corollary} {\rm (Prescription of finite structure and row minimal indices)}
  \label{corprescrfifrmi}
  Let $P(s)\in\efe[s]^{m\times n}$,  $\deg(P(s))=d$, $\rank (P(s))=r$.
Let $\alpha_1(s)\mid\cdots\mid\alpha_r(s)$ be the
  invariant factors, $e_1\leq\dots\leq e_r$ the partial multiplicities of $\infty$,
$\bc=(c_1,  \dots,  c_{n-r})$  the  column minimal indices and $ \bu=(u_1, \dots,  u_{m-r})$  the row minimal
indices of $P(s)$, where 
$u_1 \geq  \dots  \geq  u_{\eta}  > u_{\eta +1} = \dots =u_{m-r} = 0$.

Let $z$ and $x$ be integers such that $0\leq x\leq \min\{z, n-r\}$ and 
let  $\beta_1(s)\mid\dots \mid \beta_{r+x}(s) $ be monic polynomials and
$v_1 \geq  \dots  \geq  v_{\bar \eta} > v_{\bar \eta+1}=\dots= v_{m+z-r-x}=0$ 
  be non negative integers.
  Let
  \begin{equation}\label{eqbigtildeBfifrmi}
\tilde B=\sum_{i=1}^{r+x}\deg(\beta_i)-\sum_{i=1}^{r}\deg(\alpha_i)-\sum_{i=1}^{x}c_{i}+\sum_{i=1}^{m+z-r-x}v_i-\sum_{i=1}^{m-r}u_i-xd
.\end{equation}
  There exists $W(s)\in \efe[s]^{z\times n}$
 such that $\deg(W(s))\leq d$, $\rank \left(\begin{bmatrix}P(s)\\W(s)\end{bmatrix}\right) =r+x$
 and $\begin{bmatrix}P(s)\\W(s)\end{bmatrix}$  has
 $\beta_1(s)\mid\dots \mid \beta_{r+x}(s) $  as invariant factors
 and
$\bv=(v_1,\dots, v_{m+z-r-x})$ as row minimal indices
if and only if (\ref{eqetapol}),
(\ref{eqinterfifpol}),
\begin{equation}\label{eqtildeBleqfif}
\tilde B\leq \sum_{i=r+x-z+1}^{r}e_i,
\end{equation}
\begin{equation}\label{eqvustildeBfif}
 \sum_{i=1}^{m+z-r-x}v_i-\sum_{i=1}^{m-r} u_i\geq \max\{0, \tilde B\}+\sum_{i=1}^{r+x}\deg(\lcm(\alpha_{i-x}, \beta_i))-\sum_{i=1}^{r+x}\deg(\beta_i),
\end{equation}
\begin{equation}\label{eqcmimajpoltildeBfif}  (c_1, \dots, c_x) \prec   \tilde  \ba,\end{equation}
\begin{equation}\label{eqrmimajpoltildeBfif} \bv \prec'  (\bu , \tilde \bb),\end{equation}
where $\tilde  \ba = (\tilde a_1, \dots, \tilde a_x )$ and $\tilde \bb = (\tilde b_1, \dots,  \tilde b_{z-x})$ 
are 
\begin{equation}\label{eqtildeafifrmi}
\begin{array}{rl}\tilde  a_1=& \sum_{i=1}^{m+z-r-x}v_i-\sum_{i=1}^{m-r}u_i\\&+\sum_{i=1}^{r+x}\deg(\beta_i)-\sum_{i=1}^{r+x-1}\deg(\lcm(\alpha_{i-x+1}, \beta_i))-\tilde B-d
  ,\\ \tilde  a_j=&
\sum_{i=1}^{r+x-j+1}\deg(\lcm(\alpha_{i-x+j-1}, \beta_i))\\&-\sum_{i=1}^{r+x-j}\deg(\lcm(\alpha_{i-x+j}, \beta_i))
  -d,\quad 2\leq j \leq x,
\end{array}
\end{equation}
\begin{equation}\label{eqtildebfifrmi}
\begin{array}{rl}\tilde  b_1=& \sum_{i=1}^{m+z-r-x}v_i-\sum_{i=1}^{m-r}u_i+\min\{0,e_r-\tilde B \}\\&+\sum_{i=1}^{  r+x}\deg(\beta_i)-\sum_{i=1}^{  r+x}\deg(\lcm(\alpha_{i-x-1}, \beta_i)),
  \\ \tilde b_j=&\min\{0,\sum_{i=r-j+1}^{r}e_i-\tilde B\}-\min\{0,\sum_{i=r-j+2}^{r}e_i-\tilde B \}\\&
  +\sum_{i=1}^{  r+x}\deg(\lcm(\alpha_{i-x-j+1}, \beta_i))-\sum_{i=1}^{  r+x}\deg(\lcm(\alpha_{i-x-j}, \beta_i)),\\&
  \hfill 2\leq j \leq z-x.
\end{array}
\end{equation}
\end{corollary}
{\bf Proof.} 
The result is a consequence of Theorem \ref{theoprescrfifrmicmi}. The proof follows the scheme of Corollary \ref{corprescrdeirmi}  by exchanging the roles of the invariant factors and the infinite elementary divisors. 
\hfill $\Box$

\subsection{Prescription of singular structures}
In this last subsection we prescribe the row and column minimal indices, then only the row minimal indices, and finally only the column minimal indices. 
The sufficiency parts of Theorem \ref{theoprescrsing} and Corollary \ref{corprescrrmi} require the field to be algebraically closed, while 
Corollary \ref{propprescrcmi}
is valid over arbitrary fields.

  \begin{theorem} {\rm (Prescription of  singular structure)}
  \label{theoprescrsing}
  Let $\efe$ be an algebraically closed field.
  Let $P(s)\in\efe[s]^{m\times n}$,  $\deg(P(s))=d$, $\rank (P(s))=r$.
Let $\phi_1(s,t)\mid\cdots\mid\phi_r(s,t)$ be the
homogeneous invariant factors,
$\bc=(c_1,  \dots,  c_{n-r})$  the
column minimal indices, and $ \bu=(u_1, \dots,  u_{m-r})$  the row minimal
indices of $P(s)$, where 
$u_1 \geq  \dots  \geq  u_{\eta}  > u_{\eta +1} = \dots =u_{m-r} = 0$. 

Let $z$ and $x$ be integers such that $0\leq x\leq \min\{z, n-r\}$ and let $d_1 \geq  \dots   \geq d_{n-r-x} \geq 0$ and $v_1 \geq  \dots  \geq  v_{\bar \eta} > v_{\bar \eta+1}=\dots= v_{m+z-r-x}=0$ be non negative integers.
  Let
  \begin{equation}\label{eqbigEsing}
E=\sum_{i=1}^{n-r-x}d_i-\sum_{i=1}^{n-r}c_i+\sum_{i=1}^{m+z-r-x}v_i-\sum_{i=1}^{m-r}u_i-xd.\end{equation}
There exists $W(s)\in \efe[s]^{z\times n}$
such that $\deg(W(s))\leq d$, $\rank \left(\begin{bmatrix}P(s)\\W(s)\end{bmatrix}\right) =r+x$, $\begin{bmatrix}P(s)\\W(s)\end{bmatrix}$  has 
$\bv=(v_1,\dots, v_{m+z-r-x})$ as row minimal indices and $\bd=(d_1,\dots, d_{n-r-x})$ as  column minimal indices
if and only if (\ref{eqetapol}),
\begin{equation}\label{eqEleq}
E\leq \sum_{i=r+x-z+1}^{r}\deg(\phi_{i}),
\end{equation}
\begin{equation}\label{eqvusE}
 \sum_{i=1}^{m+z-r-x}v_i-\sum_{i=1}^{m-r} u_i\geq \max\{0, E\},
\end{equation}
\begin{equation}\label{eqcmimajpoltildeE}  \bc \prec'  (\bd , \tilde \ba),\end{equation}
\begin{equation}\label{eqrmimajpoltildeE} \bv \prec'  (\bu , \tilde \bb),\end{equation}
where $\tilde \ba = (\tilde a_1, \dots, \tilde a_x )$ and $\tilde \bb = (\tilde b_1, \dots,  \tilde b_{z-x})$ 
are 
\begin{equation}\label{eqasing}
\begin{array}{rl}\tilde a_1=&
{
  \sum_{i=1}^{m+z-r-x}v_i-\sum_{i=1}^{m-r}u_i-E-d,}\\ \tilde a_j=&-d,\quad 2\leq j \leq x,
\end{array}
\end{equation}
\begin{equation}\label{eqbsing}
\begin{array}{rl}\tilde b_1=& \sum_{i=1}^{m+z-r-x}v_i-\sum_{i=1}^{m-r}u_i+\min\{0,\deg(\phi_{r})-E \},
  \\ \tilde b_j=&\min\{0,\sum_{i=r-j+1}^{r}\deg(\phi_{i})-E \}-\min\{0,\sum_{i=r-j+2}^{r}\deg(\phi_{i})-E \}
  ,\\&\hfill 2\leq j \leq z-x.
\end{array}
\end{equation}
\end{theorem}
\begin{remark}\label{remac3}
The necessity of the conditions, and the sufficiency when $E\leq 0$ hold for arbitrary fields.
  \end{remark}
    {\bf Proof of Theorem \ref{theoprescrsing}}
The proof is analogous to that of Theorem \ref{theoprescrdeirmicmi} by exchanging the roles of the invariant factors and the homogeneous invariant factors. 
\hfill $\Box$

\medskip

Example \ref{exalgclosed} shows that, in general,  conditions (\ref{eqetapol}) and (\ref{eqEleq})-(\ref{eqrmimajpoltildeE})
  are not sufficient if the field is not algebraically closed and $E>0$.
  
\begin{corollary}{\rm (Prescription of row minimal indices)}
  \label{corprescrrmi}
  Let $\efe$ be an algebraically closed field.
Let $P(s)\in\efe[s]^{m\times n}$,  $\deg(P(s))=d$, $\rank (P(s))=r$.
  Let $\phi_1(s,t)\mid\cdots\mid\phi_r(s,t)$ be the homogeneous invariant factors, 
  $\bc=( c_1,\dots, c_{n-r})$  the  column minimal indices and
  $\bu=(u_1, \dots, u_{m-r})$ the  row minimal indices of $P(s)$,
  where $u_1 \geq  \dots  \geq  u_{\eta}  > u_{\eta +1} = \dots =u_{m-r} = 0$.

Let $z$ and $x$ be integers such that $0\leq x\leq \min\{z, n-r\}$ and 
let  $v_1 \geq  \dots  \geq  v_{\bar \eta}  > v_{\bar \eta +1} = \dots =v_{m+z-r-x} = 0$
be non negative integers.
Let
 $$\hat E=-\sum_{i=1}^{x}c_i+\sum_{i=1}^{m+z-r-x}v_i-\sum_{i=1}^{m-r}u_i-xd.$$ 
There exists $W(s)\in \efe[s]^{z\times n}$
such that
$\deg(W(s))\leq d$, $\rank \left(\begin{bmatrix}P(s)\\W(s)\end{bmatrix}\right) =r+x$
and $\begin{bmatrix}P(s)\\W(s)\end{bmatrix}$  has 
$\bv=(v_1,\dots, v_{m+z-r-x})$ as row minimal indices
if and only if
(\ref{eqetapol}),
\begin{equation}\label{eqEleqrm}
\hat E\leq \sum_{i=r+x-z+1}^{r}\deg(\phi_{i}),
\end{equation}
\begin{equation}\label{eqvusrm}
 \sum_{i=1}^{m+z-r-x}v_i-\sum_{i=1}^{m-r} u_i\geq \max\{0, \hat E\},
\end{equation}
\begin{equation}\label{eqrmimajpolhatsrm} \bv \prec'  (\bu , \hat \bb),\end{equation}
where $\hat \bb=(\hat b_1, \dots, \hat b_{z-x})$ is 
$$\begin{array}{rl}\hat b_1=& \sum_{i=1}^{m+z-r-x}v_i-\sum_{i=1}^{m-r}u_i+\min\{0,\deg(\phi_{r})-\hat E \},
  \\ \hat b_j=&\min\{0,\sum_{i=r-j+1}^{r}\deg(\phi_{i})-\hat E \}-\min\{0,\sum_{i=r-j+2}^{r}\deg(\phi_{i})-\hat E \},\\&
  \hfill 2\leq j \leq z-x.
\end{array}
$$
\end{corollary}
\begin{remark}\label{remac4}
The necessity of the conditions, and the sufficiency when $\hat E\leq 0$ hold for arbitrary fields.
  \end{remark}

{\bf Proof of Corollary \ref{corprescrrmi}.}
Let us assume that 
there exists a polynomial matrix  $W(s)\in \efe[s]^{z\times n}$
such that
$\deg(W(s))\leq d$, $\rank\begin{bmatrix}P(s)\\W(s)\end{bmatrix}=r+x$
and $\begin{bmatrix}P(s)\\W(s)\end{bmatrix}$  has 
$\bv=(v_1,\dots, v_{m+z-r-x})$ as row minimal indices.
 Let $\bd=(d_1,\dots, d_{n-r-x})$ be its column minimal indices.  
By Theorem \ref{theoprescrsing}, (\ref{eqetapol}) and (\ref{eqEleq})-(\ref{eqrmimajpoltildeE}) hold,
where $E$, $\tilde \ba$ and $\tilde \bb$ are defined in (\ref{eqbigEsing}), (\ref{eqasing}) and (\ref{eqbsing}), respectively.
From (\ref{eqcmimajpoltildeE}), we have $d_i\geq c_{i+x}$, $1\leq i\leq n-r-x$; hence $E\geq \hat E$.
Then from (\ref{eqEleq}) and (\ref{eqvusE}) we obtain (\ref{eqEleqrm}) and (\ref{eqvusrm}), respectively.
Moreover,
$$\begin{array}{rl}
  \sum_{i=1}^{j}\tilde b_i=&\sum_{i=1}^{m+z-r-x}v_i-\sum_{i=1}^{m-r}u_i+\min\{0, \sum_{i=r-j+1}^{r}\deg(\phi_{i})-E\}\\\leq&
  \sum_{i=1}^{m+z-r-x}v_i-\sum_{i=1}^{m-r}u_i+\min\{0, \sum_{i=r-j+1}^{r}\deg(\phi_{i})-\hat E\}\\=&\sum_{i=1}^{j}\hat  b_i
\quad 1\leq j \leq z-x.\end{array}$$
And, from (\ref{eqEleq}) and (\ref{eqEleqrm}), we derive
$\sum_{i=1}^{z-x}\tilde b_i=\sum_{i=1}^{m+z-r-x}v_i-\sum_{i=1}^{m-r}u_i=\sum_{i=1}^{z-x}\hat b_i.$
Therefore $\tilde \bb \prec \hat \bb$, and from (\ref{eqrmimajpoltildeE}) and Remark \ref{gmajcp} we obtain (\ref{eqrmimajpolhatsrm}).

Conversely, assume that  (\ref{eqetapol}) and (\ref{eqEleqrm})-(\ref{eqrmimajpolhatsrm}) hold.
Define
$d_i=c_{i+x}$, $1\leq i \leq n-r-x$, 
$\bd=(d_1, \dots, d_{n-r-x})$,
and let $E$, $\tilde \ba$ and $\tilde \bb$ be as in (\ref{eqbigEsing}), (\ref{eqasing}) and (\ref{eqbsing}), respectively.
Then
$$E=\hat E, \quad \tilde a_1=\sum_{i=1}^{x}c_i+(x-1)d, \quad \tilde \bb=\hat \bb.$$
Conditions (\ref{eqEleqrm}), (\ref{eqvusrm}), and (\ref{eqrmimajpolhatsrm}) are respectively equivalent to
(\ref{eqEleq}), (\ref{eqvusE}),  and (\ref{eqrmimajpoltildeE}).
Finally, $\bc=\bd\cup (c_1, \dots, c_x)$.
By Lemma \ref{lemmacup},
\begin{equation}\label{eqcx}
\bc\prec'(\bd, (c_1, \dots, c_x)).\end{equation}
As $\sum_{i=1}^x\tilde a_i=\sum_{i=1}^xc_i$ and
 $\sum_{i=1}^j c_i \leq \sum_{i=1}^x c_i \leq \sum_{i=1}^x c_i+(x-j)d=\sum_{i=1}^j
\tilde a_i$, $1\leq j \leq x$,  
we obtain $(c_1, \dots, c_x)\prec \tilde \ba$ and, from (\ref{eqcx}) and Remark \ref{gmajcp} we derive (\ref{eqcmimajpoltildeE}).
By Theorem \ref{theoprescrsing},
the sufficiency of the conditions (\ref{eqetapol}) and (\ref{eqEleqrm})-(\ref{eqrmimajpolhatsrm})  follows.
\hfill $\Box$

\medskip

Example \ref{exalgclosed}
shows that, in general,  conditions (\ref{eqetapol}) and (\ref{eqEleqrm})-(\ref{eqrmimajpolhatsrm})
  are not sufficient if the field is not algebraically closed and $\hat E>0$.

The following result is obtained as a consequence of Corollary \ref{corprescrdeicmi}.

\begin{corollary} {\rm (Prescription of column minimal indices)}
  \label{propprescrcmi}
Let $P(s)\in\efe[s]^{m\times n}$,
$\deg(P(s))=d$, $\rank (P(s))=r$.
  Let   
$\bc=(c_1,  \dots,  c_{n-r})$  be the
  column minimal indices of $P(s)$.
  
  Let $z$ and $x$ be integers such that $0\leq x\leq \min\{z, n-r\}$ and let 
  $d_1 \geq  \dots  \geq  d_{n-r-x}\geq 0$ 
  be non negative integers.
  There exists a polynomial matrix  $W(s)\in \efe[s]^{z\times n}$ such that  $\deg(W(s))\leq d$, $\rank\left(\begin{bmatrix}P(s)\\W(s)\end{bmatrix}\right)=r+x$ 
and $\begin{bmatrix}P(s)\\W(s)\end{bmatrix}$ has $\bd=(d_1,\dots, d_{n-r-x})$ as column minimal indices 
if and only if
\begin{equation}\label{eqdcxd}
\sum_{i=1}^{n-r-x}d_i-\sum_{i=1}^{n-r}c_i\leq xd,
\end{equation}
\begin{equation}\label{eqcmimajpolhatseq}  \bc \prec'  (\bd , \hat \ba),\end{equation}
where $\hat \ba=(\hat a_1, \dots, \hat a_{x})$ is 
$$\begin{array}{rl}
  \hat a_1 =& \sum_{i=1}^{n-r}c_i-\sum_{i=1}^{n-r-x}d_i+(x-1)d,
  \\ \hat a_j=&-d
  ,\quad 2\leq j \leq x.
\end{array}
$$
\end{corollary}
{\bf Proof.}
Let $e_1\leq\dots\leq e_r$ be the  partial multiplicities of $\infty$ of $P(s)$.

Let us assume that there exists a polynomial matrix  $W(s)\in \efe[s]^{z\times n}$ such that  $\deg(W(s))\leq d$, $\rank\left(\begin{bmatrix}P(s)\\W(s)\end{bmatrix}\right)=r+x$ 
and $\begin{bmatrix}P(s)\\W(s)\end{bmatrix}$ has $\bd=(d_1,\dots, d_{n-r-x})$ as column minimal indices. Let $f_1\leq \dots \leq f_{r+x}$
be its partial multiplicities of $\infty$.
By Corollary \ref{corprescrdeicmi}, (\ref{eqinteriedpol}),
(\ref{eqcde}) and (\ref{eqcmimajtildea}) hold, where $\tilde \ba$ is defined in (\ref{eqaiedcmi}).

Since
$\sum_{i=1}^{r+x-j}\max\{e_{i-x+j}, f_i\}\geq \sum_{i=1}^{r}e_i$, $0\leq j \leq x$,
 (\ref{eqcde}) implies (\ref{eqdcxd}) and
$$
\begin{array}{rl}
\sum_{i=1}^{j} \tilde a_i=&\sum_{i=1}^{n-r}c_i-\sum_{i=1}^{n-r-x}d_i+(x-j)d\\&+\sum_{i=1}^{r}e_i-\sum_{i=1}^{r+x-j}\max\{e_{i-x+j}, f_i\}\\
\leq&\sum_{i=1}^{n-r}c_i-\sum_{i=1}^{n-r-x}d_i+(x-j)d=\sum_{i=1}^{j} \hat a_i, \quad 1\leq j \leq x.
\end{array}
$$
Moreover,  from  (\ref{eqinteriedpol}) we get
$
\sum_{i=1}^{x} \tilde a_i=\sum_{i=1}^{x} \hat a_i.
$
Thus, $\tilde \ba \prec \hat \ba$,  and  (\ref{eqcmimajtildea}) and Remark \ref{gmajcp} imply (\ref{eqcmimajpolhatseq}).

Conversely, assume that  (\ref{eqdcxd}) and (\ref{eqcmimajpolhatseq}) hold.
Define
$f_i=e_{i-x}$, $1\leq i \leq r+x,$
and $\tilde \ba$ as in (\ref{eqaiedcmi}).
Then (\ref{eqinteriedpol}) holds, and
$\sum_{i=1}^{r+x-j}\max\{e_{i-x+j}, f_i\}= \sum_{i=1}^{r}e_i$, $0\leq j \leq x;$
hence (\ref{eqdcxd}) and (\ref{eqcmimajpolhatseq}) are equivalent to  (\ref{eqcde}) and (\ref{eqcmimajtildea}), respectively.
By Corollary \ref{corprescrdeicmi}, the sufficiency of the conditions (\ref{eqdcxd}) and (\ref{eqcmimajpolhatseq}) follows.
\hfill $\Box$

\begin{example}
  \label{exprifcmi} 
  Let $P(s)=
  \begin{bmatrix}s^2&-1&0\\0&0&0
  \end{bmatrix}
  \in \FF[s]^{2\times3}
  $. Then
  $d=2$,  $r=1$, $\bc=(2,0)$.
  Let $x=z=1$. 
  In this case, condition (\ref{eqdcxd}) becomes $d_1\leq 4$, and 
(\ref{eqcmimajpolhatseq})  is equivalent to $d_1\geq 2$ or $d_1=0$. 
For example, if we prescribe $d_1=2$, a possible completion is 
$Q(s)=\begin{bmatrix}s^2&-1&0\\0&0&0\\0&0&1
  \end{bmatrix}$.
  \end{example}

\section{Conclusions}

In this work the study of the row (column) completion problem  of polynomial matrices of given degree started in \cite{AmBaMaRo23} has been completed.
We  present here a solution to Problem \ref{problem} in the  cases  not studied in \cite{AmBaMaRo23}. Thus, Problem \ref{problem} has been totally solved. 
 The necessity of the conditions  is valid over arbitrary fields in all cases. Sufficiency is also valid over arbitrary fields, except in 4 cases: in the  prescription of infinite and singular structures (Theorem \ref{theoprescrdeirmicmi}), of  infinite structure and row minimal indices (Corollary \ref{corprescrdeirmi}), of singular structure (Theorem \ref{theoprescrsing}), and of row minimal indices (Corollary \ref{corprescrrmi}). In these cases the sufficiency is only valid, in general,  over algebraically closed fields, as shown in Example \ref{exalgclosed}.

  \bibliographystyle{acm} 
\bibliography{references}
\end{document}